\documentclass[12pt,a4paper]{amsart}
\usepackage[top=2.5cm,bottom=2.5cm,left=2.5cm,right=2.5cm]{geometry}%\geometry{papersize={14cm,10.5cm}}
\usepackage{cases}
\usepackage{graphicx}
\usepackage{stmaryrd}
\usepackage[all,pdf]{xy}
\usepackage{amssymb}
\usepackage{mathrsfs}
\usepackage{ulem}
\usepackage{xcolor}

\newtheorem{theorem}{Theorem}[section]
\newtheorem{corollary}[theorem]{Corollary}
\newtheorem{lemma}[theorem]{Lemma}
\newtheorem{proposition}[theorem]{Proposition}
\newtheorem{definition}[theorem]{Definition}

\numberwithin{equation}{section}

       \newcommand{\beqm}{\begin{eqnarray*}}
       \newcommand{\eqm}{\end{eqnarray*}}

\newcommand{\Cn}{\mathbb{C}^n}
\newcommand{\D}{\mathbb{D}}

\newcommand{\fr}{\frac}

\newcommand{\B}{\mathbb{B}}

\begin{document}
\title{ Absolutely Summing  Toeplitz operators on Bergman spaces in the unit ball of $\mathbb{C}^n$}

\author{Zhangjian Hu, Ermin Wang}

\address{Zhangjian Hu: Department of Mathematics, Huzhou University, Huzhou 313000, Zhejiang  China}

\email{huzj@zjhu.edu.cn}

\address{Ermin Wang: School of Mathematics and Statistics, Lingnan Normal University, Zhanjiang 524048, Guangdong  China}

\email{wem0913@sina.com}

\thanks{}



\maketitle
\begin{abstract}
 In this paper, for $p> 1 $ and $r \ge 1$ we provide a complete characterization of the positive Borel measures  $\mu$ on the unit ball $\B_n$ of $\mathbb {C}^n$  for which the induced Toeplitz operator $T_\mu$ is $r$-summing on the Bergman space $A^{p}$. We prove that the $r$-summing norm of $T_\mu: A^p\to A^p$ is equivalent to  $\|\widetilde{\mu}\|_{L^{\kappa}(d\lambda)}$, where $\kappa$ is a positive number determined by $p$ and $r$.
As some preliminary, we describe when a Carleson embedding $J_\mu: A^p \to L^q(\mu) (1\le p,  q\le 2)$ is $r$-summing,  which extends  the main result in [B. He, et al, Absolutely summing Carleson embeddings on Bergman spaces, Adv. Math., 439, 109495 (2024)].
\end{abstract}

 \section{  Introduction and Main Results}

Let $\mathbb{C}^n$ denote the $n$-dimensional complex Euclidian space. For $z=(z_1, \cdot\cdot\cdot, z_n)$ and  $w=(w_1, \cdot\cdot\cdot, w_n)$ in $\mathbb{C}^n$, we write
$\langle z,w \rangle =z_1\overline{w}_1+\cdot\cdot\cdot+z_n\overline{w}_n$ and
$|z|=\sqrt{\langle z,z \rangle}$. The open unit ball in $\mathbb{C}^n$ is the set
$${\B_n}=\{z\in \mathbb{C}^n: |z|<1\}.$$ We use $H(\B_n)$ to denote the space  of holomorphic functions in ${\B_n}$.

For $1\le p < \infty$, the  usual space  ${L}^{p}$  consists of Lebesgue  measurable functions on $\B_n$ such that
$$
\|f\|_{L^p}= \left(\int_{\B_n}|f(z)|^pdv(z)\right)^{\fr 1p}< \infty,
$$
where  $dv$ is the volume measure on $\B_n$, normalized so that $v(\B_n)=1$. The Bergman space ${A}^{p}$ is defined by
$$
A^p= H(\B_n)\cap L^p
$$
with inherited norm  written as $\|\cdot\|_{A^p}$.

 It is clear that \(A^{p}\) is a Banach space under \(\|\cdot\|_{A^{p}}\) for \(1 \leq  p < \infty\).  The  point evaluations $L_z: f\mapsto f(z)$ are bounded linear functionals on $A^p$ for each $z\in \B_n$. In particular, ${A}^{2}$ is a  reproducing kernel Hilbert space: for each $z\in \B_n$, there is a function $K_{z}\in {A}^{2}$ with $\|L_z\|=\|K_{z}\|_{A^2}$,  such that $$f(z)=\int_{\B_n}f(w)\overline{K_z(w)}dv(w)$$ for all $f\in A^{2}$.
The function $K(z, w)=\overline{K_z(w)}$ is called the Bergman kernel of  $A^2$ which has the property that $K(z, w)=\overline{K(w, z)}$.
 The orthogonal projection \(P: {L}^2  \rightarrow  A^{2}\) can be represented as
\begin{equation} \label{projection}
Pf(z) = \int_{\B_n}K(z,w) f(w)dv(w), \ \ \ z\in \B_n.
\end{equation}
 The Bloch space $\mathcal{B}$ is  defined to be the space of holomorphic functions $f$ on $\B_n$ such that
 $$
 \|f\|_\mathcal{B}=|f(0)|+\sup \left \{(1-|z|^2)\sum_{j=1}^n \left |\frac{\partial f}{\partial z_j} (z)\right |: z\in\B_n \right \}<\infty.
 $$
All the spaces $A^p$ and $\mathcal{B}$ are Banach spaces.
 For $1< p<\infty$, the dual space of ${A}^{p}$  can be identified with ${A}^{p'}$ under the pairing
$$
 \langle f, g \rangle =\int_{\B_n}f(z)\overline{g(z)}dv(z),\ \ \ f\in A^{p}, g\in A^{p'},
$$
  where $p'$ is the conjugate exponent of $p$, that is, $\frac 1p+ \frac 1{p'} =1$.
   And
 the dual space of ${A}^{1}$  can be identified with $\mathcal{B}$ under the pairing
$$
 \langle f, g \rangle = \lim_{r\to 1^-} \int_{\B_n}f(rz)\overline{g(z)}dv(z),\ \ \ f\in A^{1}, g\in \mathcal{B}.
$$
 See \cite{Zh05} for more detail.

With the Bergman kernel $K(z, w)$ we define the Bergman metric $B(z) $ on $\B_n$ as
$$
   B(z)= \frac 1{n+1} \left( \frac {\partial^2}{\partial \bar{z_j} \partial
   z_k} \log K(z, z)  \right)_{n\times n}.
$$
The induced Bergman distance  is denoted by $\beta(\cdot, \cdot)$.
 Given  $z\in \B_n $ and $r>0$, let $B(z, r)$  denote the Bergman ball  centered at $z$  with radius $r$, that is,
$$
B(z, r)=\{w\in \B_n: \beta(w,z)<r\}.
$$
For $1< p<\infty$ and $z\in \B_n$, write
$
k_{z, p}(\cdot)= {K_z(\cdot)}/{\|K_z\|_{A^p}}
$
to be  the normalized reproducing kernel for $A^{p}$. Given some positive Borel measure $\mu$ on $\B_n$ (denoted by $\mu \ge 0$), the Berezin transform of $\mu$ is defined as
$$
\widetilde{\mu}(z)=\int_{\B_n} \left| k_{z,2}(w)\right|^2 d\mu(w).
$$
 For fixed  $\delta >0$, set
$$
\widehat{\mu}_\delta(z)= \frac{\mu(B(z, \delta))}{v(B(z, \delta))}, \ \ \ z\in \B_n,
$$
to be the average of  $\mu $ on $B(z, \delta)$.

The Toeplitz operator $T_\mu$ with symbol $\mu\ge 0$ is defined as
\begin{equation}\label{toep}
T_\mu f(z) = \int_{\B_n} f(w)K(z,w)d\mu(w).
\end{equation}
if it is  densely well defined on $A^p$.\\

Here is the definition of $r$-summing operators between Banach spaces. \

\begin{definition} \label{def1}
Suppose \(1 \leq  r< \infty\) , a linear operator \(T : X \rightarrow  Y\) between Banach spaces \(X\) and \(Y\) is said to be \(r\)-summing, written as $T\in \Pi_r(X, Y)$, if there exists \(C >  0\) such that
\begin{equation}\label{Main-def}
{\left( \mathop{\sum }\limits_{{k = 1}}^{n}{\begin{Vmatrix}T{x}_{k}\end{Vmatrix}}_{Y}^{r}\right) }^{\fr 1r} \leq  C\mathop{\sup }\limits_{{{x}^{ * } \in  {B}_{{X}^{ * }}}}{\left( \mathop{\sum }\limits_{{k = 1}}^{n}{\left| {x}^{*}({x}_{k}) \right| }^{r}\right) }^{\fr 1r}
\end{equation}
for every finite sequence \({\left\{  {x}_{k}\right\}  }_{1 \leq  k \leq  n} \subset  X\), where ${X}^{ * }$ denotes the dual space of $X$, and ${B}_{{X}^{*}}$ denotes the closed unit ball $\{x^*\in {X}^{*}: \|x^*\|\le 1\}$ of  ${X}^{ * }$.
\end{definition}

Denote by $\Pi_r(X,Y)$ the set of all $r$-summing operators from $X$ to $Y$. The $r$-summing norm of $T$, denoted $\pi_r(T: X \to Y)$, is defined as the infimum over all constants $C$ satisfying inequality (\ref{Main-def}). Clearly, $\Pi_r(X,Y)$ is a linear subspace of $B(X,Y)$, the space of bounded linear operators from $X$ to $Y$, and $\pi_r$ defines a norm on $\Pi_r(X,Y)$ satisfying $$\|T\|\le \pi_r(T: X\to Y).$$  Moreover, it is well known that $\Pi_r(X,Y)$ is a Banach space under the norm $\pi_r$. In what follows, we abbreviate $\pi_r(T: X \to Y)$ as $\pi_r(T)$ whenever no confusion occurs.\\

The theory of absolutely $r$-summing operators plays an important role in the general theory of Banach spaces \cite{LP68}. The study of absolutely $1$-summing operators was initiated by Grothendieck \cite{Gr53, Gr55}, who defined them via tensor products with the space $\ell^1$. One of his celebrated results in this direction is that every bounded linear operator from $l^1$ to $l^2$ is $1$-summing-a statement now known as Grothendieck's theorem. Pietsch \cite{Pi67} later introduced the notion of absolutely $r$-summing operators between normed spaces and established many of their fundamental properties. Such operators naturally generalize Hilbert-Schmidt operators on Hilbert spaces to the setting of general Banach spaces. Indeed, for Hilbert spaces $H_1$ and $H_2$, an operator $T: H_1 \to H_2$ is Hilbert-Schmidt if and only if it is $2$-summing, and moreover $|T|_2 = \pi_2(T)$.

In the recent years, the study of $r$-summing operators on holomorphic function spaces has attracted considerable attention. In 2018, on the unit disc $\D$ Lefevre and Rodriguez-Piazza \cite{LR18} obtained a complete characterization of $r$-summing Carleson embeddings on Hardy spaces $H^p(\D)$. Building on this work, analogous characterizations were established for embeddings on the classical Bergman spaces $A^p(\mathbb{D})$ \cite{HJL24}, weighted Bergman spaces $A^p_\alpha(\mathbb{H})$ over the upper half-plane $\mathbb{H}$ \cite{CDW25}, and weighted Fock spaces $F^p_{\alpha,\omega}(\mathbb C)$ with $A_\infty$-type weights \cite{CHW24}. Additionally, absolutely summing Volterra operators and weighted composition operators on Bergman and Bloch spaces
have been investigated in \cite{FL22, JL24}. \\

Toeplitz operators are widely recognized for their significant applications in fields such as signal processing and quantum theory,  they also  play a key  role in the operator theory of holomorphic function spaces. There has been extensive research on Toeplitz operators acting on Bergman spaces from various perspectives. Luecking  \cite{Lu87} was probably the pioneer to study Toeplitz operators $T_\mu$ with measures as symbols, where he provides, among other things, a description of Schatten class Toeplitz operators $T_\mu: A^2_\alpha\to A^2_\alpha$ in terms of an $l^p$-condition involving a hyperbolic lattice of $\mathbb{D}$, where $A^2_\alpha$ is the weighted Bergman space as (\ref{weighted}). Subsequent works have further explored positive Toeplitz operators on various weighted Bergman spaces; see, for example, \cite{DHLZ25, PZ15, PRS16, WW22, Zh88, ZWH23} and the references therein.

It is worth noting that research on Toeplitz operators has predominantly focused on their boundedness, compactness, Schatten class membership, and related algebraic properties. However, the study of absolutely summing Toeplitz operators on Bergman spaces remains unexplored to date. The aim of this paper is to fill this gap by providing a comprehensive characterization of positive Toeplitz operators that are
$r$-summing on the Bergman space.\\

The purpose of the present work  is the characterize those    $\mu\ge 0$ on $\B_n$ so that $T_\mu$ is $r$-summing on the Bergman space $A^p$ provided    $1\le p,  r<\infty$. To state our main result precisely,  we  write
 $$
 d\lambda(z)=K(z,z)dv(z)
 $$
to be  the invariant measure on  $\B_n$.
 Given $1\le p, r<\infty$ we set the constant $\kappa=\kappa(p, r)$ as
 $$
  \kappa=
\left\{
  \begin{array}{ll}
   2, \ \ \ \ \textrm{if} \ 1\le p\le 2; \\
    p', \ \ \ \textrm{if} \ p\ge 2 \ \textrm{ and }\ 1\le r\le p';\\
    r, \, \ \ \ \textrm{if} \ p\ge 2 \ \textrm{ and }\ p'\le r\le p;\\
    p, \, \ \ \ \textrm{if} \ p\ge 2 \ \textrm{ and }\ p\le r<\infty.
  \end{array}
\right.
$$
Two quantities $A$ and $B$ are said to be equivalent, written as  $A\simeq B$, if  $C^{-1}A \le B \le CA$ for some $C>0$.

\begin{theorem}\label{THM2}
 Let $\mu$ be a positive Borel measure on \(\B_n\). For $1< p<\infty$ and $1\le
  r<\infty$,
the Toeplitz operator   \(T_\mu: A^p \rightarrow  A^p\) is \(r\)-summing if and only if \(\widetilde{\mu} \in  {L}^{\kappa}(d\lambda)\). Furthermore,
$$
\pi_r(T_\mu : A^p \rightarrow  A^p) \simeq  \|\widetilde{\mu}\|_{L^{\kappa}(d\lambda)}.
$$
The  statement above  remains true  if $\widetilde{\mu}$ is replaced by  $\widehat{\mu}_\delta$ with some (or any)  $\delta>0$.
\end{theorem}

In \cite{HW25}, analogous problems were studied in the setting of Fock spaces. It should be noted that the methods employed in the present paper differ substantially from those in \cite{HW25}. Among these differences, two points are particularly crucial. First, the Fock spaces considered in \cite{HW25} satisfy the nesting property: $F^p_\varphi \subset F^q_\varphi$ for $0 < p \le q$. Consequently, the embedding $\mathrm{Id}: F^1_\varphi \to F^2_\varphi$ is bounded whenever $p \le q$, and Grothendieck’s theorem immediately implies that $\mathrm{Id}: F^1_\varphi \to F^2_\varphi$ is $1$-summing. Such a nesting property does not hold in the Bergman space setting. Second, the characterization for the $r$-summing Toeplitz operator $T_\mu$ on $F^1_\varphi$ coincides with that on $F^2_\varphi$. This uniformity allows one to use complex interpolation to treat $r$-summing Toeplitz operators on $F^p_\varphi$ for $1 \le p \le 2$. By contrast, the boundedness characterization for $T_\mu$ on the Bergman space $A^1$ differs from that on $A^p$ ($p > 1$). This suggests that whether $T_\mu$ is $r$-summing on $A^1$ may differ from its behavior on $A^p$. Consequently, complex interpolation cannot be applied  to the Bergman space case.\\

As one preliminary to prove Theorem \ref{THM2},  we have the following theorem on the $r$-summing Carleson embedding.

\begin{theorem}\label{THM1}
Suppose $1\le p, q\le 2$ and $ s=\frac{2p}{2p-2q+pq}$.  Let $\mu$ be a positive Borel measure on \(\B_n\). Then for  \(r \geq  1\), the embedding 
$$
J_\mu:  A^p \rightarrow  L^q(\mu)), \ \ f\mapsto f
$$
is \(r\)-summing if and only if
  \(\widetilde{\mu}(\cdot)  K(\cdot, \cdot)^{\frac{q}{2}} \in  {L}^s\).
Furthermore.
\begin{equation}\label{norm121}
\pi_r(J_\mu: A^p \rightarrow  A^p) \simeq  \left\|\widetilde{\mu} (\cdot)K(\cdot, \cdot)^{\frac{q}{2}}\right\|_{L^s}^{\fr1q}. 
\end{equation}
The  statement above  remains true  if $\widetilde{\mu}$ is replaced by  $\widehat{\mu}_\delta$ with some (or any)  $\delta>0$.
\end{theorem}

Theorem \ref{THM1} extends the main result of \cite{HJL24} for $1<p\le 2$ to a larger  range $1\le p\le 2$  and to arbitrary dimension $n$.  Moreover, by avoiding the heavy machinery invoked in \cite{HJL24}, our proof is both simpler and more direct, and we expect it to work equally well on bounded symmetric domains and on bounded strongly pseudoconvex domains; see Section \ref{remark}.\\

The paper is organized as follows. In Section 2, we present some preliminaries that will be used throughout. In Section 3, we are going to prove Theorem \ref{THM1}.   Section 4 is devoted to the proof of Theorem \ref{THM2}.  Starting from Pietsch factorization theorem, we exploit the relationship between absolutely summing operators and order bounded operators. Combining this with Grothendieck's theorem and Pietsch domination theorem yields the desired conclusion. In Section 5, we would like to give some remarks concerning certain further generalization.

 Throughout this paper, we use $C$ to denote positive constants whose value may change from line to line, but do not depend  on functions being considered. For two quantities $A$ and $B$ we write  $A\lesssim B$ if there exists some $C$ such that
    $A\leq CB$.

 \section{Preliminaries}
In this section, we collect  some notations and preliminary results that  will be used throughout the paper.

Let $\delta>0$ and let $\{a_k\}_{k=1}^{\infty}$ be some sequence in  $\B_n$.  We call $\{a_k\}_{k=1}^{\infty}$ a  $\delta$-lattice in the Bergman metric if  $\left\{B(a_k,\delta)\right\}_{k}$ covers $\B_n$ and  there is some constant $c>0$ so that  $\left\{B\left(a_k,c\delta \right)\right\}_{k}$ are pairwise disjoint.
The existence of $\delta$-lattices can be seen from Theorem 2.23 in \cite{Zh05}.     Given a $\delta$-lattice $\{a_k\}_{k=1}^{\infty}$ we have an integer $N$ so that
$$
   1\le \sum_{k=1}^\infty \chi_{B(a_k,\delta)}(z)\le N
$$
for $ z\in {\B_n}$, where $\chi_{E}$ is the characteristic function of  $E\subset \B_n$.

 Let $\{r_k\}_{k=1}^\infty$  be  the sequence of Rademacher functions on  $[0, 1]$, defined by
$$ r_{1}(t)=
\left\{
  \begin{array}{ll}
    1, \ \ \ \ \ \textrm{if} \ 0\leq t-[t]<\frac{1}{2}; \\
    -1, \ \ \ \textrm{if} \ \frac{1}{2}\leq t-[t]<1,
  \end{array}
\right.
$$
and
\begin{equation}\label{rademacher}
r_{k+1}(t)=r_{1}(2^{k}t) \ \textrm{ for} \ \ k=1,2, \cdots.
 \end{equation}
 The Khinchine's inequality is closely related to the  Rademacher functions, which plays an important role in our study. For $0<l<\infty$,  there exists some positive constants $C_{1}$  and $C_{2}$ depending only on $l$ such that
$$
C_{1}
\left(\sum_{k=1}^{m}|b_{k}|^{2}\right)^{\frac{l}{2}}\leq\int_{0}^{1}\left|\sum_{k=1}^{m}b_{k}r_{k}(t)\right|^{l}dt \leq C_{2}
\left(\sum_{k=1}^{m}|b_{k}|^{2}\right)^{\frac{l}{2}}
$$
for all $m\geq 1$ and complex numbers $b_{1}, b_{2}, \ldots, b_{m}$. See \cite{DJT95} for details.

Recall that for $1\le p<\infty$, the space $l^p$ is given by
$$
l^p=\left\{x=\{x_n\}_{n\ge 1}\subset \mathbb{C}^n: \|x\|_{l^p}=\left(\sum_{n=1}^{\infty}|x_n|^p\right)^{\fr 1p}<\infty\right\}.
$$

%Let $\beta$ be the Bergman metric is defined by
% $$
%\beta(z,w)=\fr12{\rm log}\frac{1+\rho(z,w)}{1-\rho (z,w)},
%$$
%where $\rho(z,w)$ is the pseudo-hyperbolic metric on $\D$,  defined by
%$$
%\rho (z,w)=\left|\frac{z-w}{1-z\overline{w}}\right|, \ \ \ z,w \in \D.
%$$

The following two lemmas can be found in  \cite{Zh05, Zh07}.

 \begin{lemma}\label{lem3}
Let $\{a_k\}_{k=1}^\infty$ be an $\delta$-lattice with $\delta>0$,  then for  $1\le p< \infty$ and $\{c_k\}\in l^p$, the function  $f$ defined by
$$
f(z)=\sum^{\infty}_{k=1}c_kk_{a_k, p}(z), \ z\in \mathbb{B}_n
$$
belongs to $A^p$ with $\|f \|_{A^p}\lesssim \|\{c_k\} \|_{l^p}$.
\end{lemma}

\begin{lemma}\label{lem1}
Suppose $f\in A^p$, where $p>0$. Then
 \begin{equation}\label{formula1}
|f(z) |\le \|f\|_{A^p}\cdot K(z,z)^{\fr 1 p}
\end{equation}
for all $z\in \mathbb{B}_n$.
\end{lemma}

It is well know that, for $f\in H(\B_n)$, $f$ has a homogeneous expansion as
$$
     f(z)= \sum_{j=0}^\infty f_j(z). $$
Based on this expansion, for $N$ positive  we can define an invertible operator
$$
    R^{0, N}: H(\B_n)\to H(\B_n)
$$
 as
$$
    R^{0, N}f(z)= \sum_{j=0}^\infty \frac{\Gamma(n+1)\Gamma(n+1+j+N)}{\Gamma(n+1+N)\Gamma(n+1+j)} f_j(z).
$$
The inverse of $R^{0, N}$, denoted by $R_{0, N}$  is given by
$$
    R_{0, N}f(z)=\sum_{j=0}^\infty  \frac{\Gamma(n+1+N)\Gamma(n+1+j)}{\Gamma(n+1)\Gamma(n+1+j+N)} f_j(z).
$$

The following lemma comes from Lemma 1.14 in \cite{Zh05}.
\begin{lemma}\label{R-op}
For $N>0$, the operator $ R^{0, N}$ is the unique continuous linear operator on $H(\B_n)$ satisfying
$$
  R^{0, N} \left(\frac 1{(1-\langle z, w\rangle)^{n+1}} \right)= \frac {1}{(1-\langle z, w\rangle)^{n+1+N}},
$$
for all $w\in \B_n$. Similarly, the operator $ R_{0, N}$ is the unique continuous linear operator on $H(\B_n)$ satisfying
$$
R_{0, N} \left(\frac 1  {(1-\langle z, w\rangle)^{n+1+N}}\right)= \frac {1}{(1-\langle z, w\rangle)^{n+1}}.
$$
 \end{lemma}

 For $0<p,q<\infty$, we call a positive Borel measure $\mu$ on $\B_n$ is a $(p, q)$-Carleson measure if there exists a positive constant $C$ such that
  $$
  \int_{\B_n} |f(z)|^qd\mu(z)\le C\|f\|_{L^p}^q
  $$
for all $f\in A^p$.  Luecking's characterization tells us  that the $(p, q)$-Carleson measure   depends only on the ratio $s=\frac {p}{q}$, so we often call it $s$-Carleson measure.  The $s$-Carleson measures were first studied by Hastings \cite{Ha75}, and further pursued by many authors. The details can be found in \cite{Lu83, Lu85, PZ15, ZZ08, Zh07} and the references therein. In the following lemma, we list some relevant results that we need.

\begin{lemma}\label{lem2.3}
Suppose $\mu\ge 0$, then  for $1\le p<\infty$,  the following statements  are equivalent:

 {\rm (1)} $\widetilde{\mu}\in L^p(d\lambda)$;

 {\rm (2)} $\widehat{\mu}_\delta\in L^p(d\lambda)$ for some (or any) $\delta>0$;

  {\rm (3)} $\{\widehat{\mu}_\delta(a_k)\}\in l^p$ for some (or any) $\delta$-lattice $\{a_k\}$.
Moreover,
$$
\|\widetilde{\mu}\|_{L^p(d\lambda)}\simeq\|\widehat{\mu}_\delta\|_{L^p(d\lambda)}\simeq\|\{\widehat{\mu}_\delta(a_k)\}\|_{l^p}.
$$
\end{lemma}

Now we recall some useful  properties of $r$-summing operators, which can be found in \cite{DJT95}.

{\rm (1)} \textbf{Ideal property of $p$-summing operators} {\rm[p.37]}: For any $r \geq  1\) , the class ${\Pi }_{r}\left( {X,Y}\right)$ forms an operator ideal between Banach spaces: for any $T \in  {\Pi }_{r}\left( {X,Y}\right)$ , and for any two Banach spaces ${X}_{0},{Y}_{0}$ such that both $S : {X}_{0} \rightarrow  X$ and $U : Y \rightarrow  {Y}_{0}$ are linear bounded operators, we have ${UTS} \in  {\Pi }_{r}\left( {{X}_{0},{Y}_{0}}\right)$ with
\begin{equation}\label{ideal-a}
{\pi }_{r}\left( {UTS}\right)  \leq  \|U\| {\pi }_{r}\left( T\right) \| S\|.
\end{equation}

{\rm (2)} \textbf{Inclusion relations} {\rm[p.39]}: If  $1 \le p<q<\infty$, then ${\Pi }_{p}\left( {X,Y}\right)\subset{\Pi }_{q}\left( {X,Y}\right)$. Moreover, for $T \in  {\Pi }_{p}\left( {X,Y}\right)$, we have
%\begin{equation}\label{inclusion}
$$
{\pi }_{q}(T)\le {\pi }_{p}(T) \ \ \textrm{for } \ \ q\ge p.
$$
%\end{equation}

{\rm (3)} \textbf{Cotype property} We call  a Banach space $X$ has cotype $q$ if there is a constant $\kappa > 0$ such that no matter how we select finitely many vectors $x_1, x_2, \cdots, x_m$  from $X$,
$$
    \left(\sum_{j=1}^m \|x_j\|^q\right)^{\frac 1 q}\le \kappa \left( \int_0^1 \left\| \sum_{j=1}^m r_j(t) x_j\right\|^2dt\right)^{\frac 12},
$$
where $r_j$ is as in (\ref{rademacher}).

It is well known that,  for $1\le p<\infty$, the Lebesgue  space
$
L^p(\Omega, v)\ \ \textrm{ has cotype}\ \ \max\{p,2\},
$
see  \cite[p.219]{DJT95}.   For Banach spaces  $X$ and $Y$,
if $X$ and $Y$ both have cotype 2,   \cite[Corollary 11.16]{DJT95} tells us that
\begin{equation} \label{cotype}
{\Pi }_{r}\left( {X,Y}\right)={\Pi}_{1}\left( {X,Y}\right).
\end{equation}
with the equivalence $\pi_r(T)\simeq\pi_1(T)$ for $1\le r<\infty$.

The following result, known as  Pietsch domination theorem, plays a key role in our study.

\begin{lemma}\label{lem-Pietsch}
 Suppose \(1\le r<\infty\), $T: X\to Y $ is a Banach space linear operator. Then $T$
 is $r$-summing if and only if there exists a constant $C$ and a  regular probability measure $\sigma$ on the closed unit ball ${B}_{X^\ast}$ of  ${X^\ast} $ such that for each $f\in X$, there holds
 \begin{equation}\label{ineq-1}
\|Tf\|_{Y}\leq  C\cdot\left(\int_{B_{X^\ast}}{\left|\left\langle g, f \right\rangle\right| }^r{d\sigma(g)}\right)^{\fr 1 r}.
 \end{equation}
 In such a case, $\pi_r(T)$ is the least of all the constants $C$ for which such a measure exists.
\end{lemma}

We shall use several times the following well known fact about summing operators, see \cite{DJT95} or \cite{LR18}.
\begin{lemma}\label{lem8}
  Let \(X,Y\) be Banach spaces and \(r \geq  1\). Then a bounded operator $T: X\to Y $ is $r$-summing if and only if there exists $C > 0$ such that  for any measurable space $(\Omega, \sum, \nu)$ and  any continuous mapping  $F: \Omega\to X$, there holds
\[
\left (\int_{\Omega }\|T \circ  F\|_{Y}^{r}{d\nu}\right)^{\frac 1 r} \leq  C  \mathop{\sup }\limits_{{\xi  \in  {B}_{{X}^{ * }}}}\left ( \int_{\Omega }{\left| \xi  \circ  F\right| }^{r}{d\nu}\right)^{\frac 1 r}.
\]
Moreover, the best $C$ is $ {\pi }_{r}{\left( T\right) } $. 
\end{lemma}

The order bounded operators are closely related to absolutely summing operators.  Recall that, for $p\ge 1$, a Banach space linear operator $T: X\to {L}^{p}(\mu)$ is order bounded if $T(B_X)$ is an order bounded subset of $L^p(\mu)$. That is, there exists a non-negative function $h\in L^p(\mu)$ such that $|f|\le h$ $\mu$-almost everywhere for each $f\in T(B_X)$. The following lemma  is a consequence of Propositions 5.5 and 5.18 in \cite{DJT95}.

\begin{lemma}\label{lem66}
Let $X$ be a Banach space, $p \geq  1$, and let $(\Omega ,\Sigma ,m)$ be a measure space. If the operator $T : X \rightarrow  {L}^{p}\left( {\Omega ,m}\right)$ is order bounded, then it is p-summing with
$$
{\pi }_{p}\left( T\right)  \leq  {\left\|\mathop{\sup }\limits_{f \in  B_X}|Tf| \right\|}_{{L}^{p}\left( {\Omega ,m}\right) }.
$$
\end{lemma}

\section{Proof of Theorem \ref{THM1}}

We begin with the introduction of some auxiliary function spaces to be employed later.
Given $\alpha>-1,$ we define a positive Borel measure $dv_\alpha$ on $\B_n$ as follows:
$$
dv_\alpha(z)=c_\alpha(1-|z|^2)^\alpha dv(z),
$$
where  $$c_\alpha=\frac{\Gamma(n+\alpha+1)}{n!\ \Gamma(\alpha+1)}$$ is a normalizing constant so that $dv_\alpha$ is a probability measure on ${\B_n}$.
Let $\mu$ be a positive Borel measure on $\B_n$ and \(1\le p < \infty\), the space \({L}^{p}(\mu)\) consists of all Lebesgue measurable functions \(f\) on $\B_n$ for which
\[
\|f\|_{L^p(\mu)} = \left(\int_{\B_n}| f(z)|^{p}d\mu(z) \right)^{\frac{1}{p}} < \infty.
\]
When $d\mu=dv_\alpha$, we use $L^p_\alpha$ instead of $L^p(dv_\alpha)$. For $1\le p < \infty$, the  Bergman space ${A}^{p}_\alpha$ is defined by
\begin{equation}\label{weighted}
A^p_\alpha= L^p_\alpha \cap H(\mathbb{B}_n)
\end{equation}
with inherited norm written  as $\|\cdot\|_{A^p_\alpha}$.
The reproducing kernel of $A^2_\alpha$ is given by
$$
K_\alpha(z, w)=\frac {1}{(1-\langle z, w\rangle)^{n+1+\alpha}},\ \ z,w\in \B_n.
$$
When $\alpha$=0, we write $K_\alpha(z, w)$ as $K(z, w)$-the reproducing kernel of $A^2$.
Similar to (\ref{formula1}) we have 
\begin{equation}\label{weighted-a}
|f(z) |\le  \|f\|_{A^p_\alpha}\cdot K_\alpha(z,z)^{\fr 1 p}
\end{equation}
for all $f\in A^p_\alpha$ and $z\in \mathbb{B}_n$ . See \cite{Zh05}.

\begin{lemma}\label{small-a}
Let $1<p\le 2$ and let $\mu$ be a positive Borel measure on $\B_n$. The the following statements are equivalent. 

(1) The embedding operator $J_\mu: A^p \to L^2(\mu)$ is 2-summing.

(2) $\int_{\B_n} \|K(\cdot, w)\|_{L^2(\mu)}^{p'} dv(w)<\infty$. \\
Furthemore, 
\begin{equation}\label{small-b}
    \pi_2\left ({J_\mu:A^p\to L^2(\mu)}\right) \simeq  \left ( \int_{\B_n} \|K(\cdot, w)\|_{L^2(\mu)}^{p'} dv(w) \right)^{\frac 1{p'}}.
 \end{equation}    
\end{lemma}

\begin{proof}
   Suppose $J_\mu: A^p \to L^2(\mu)$ is 2-summing first. Since $p'>  2$ we know $J_\mu \in \Pi_{p'}(A^p, L^2(\mu))$ with the estimate
$$
      \pi_{p'}(J_\mu) \le \pi_2(J_\mu).
$$
For $w\in \B_n$,  take 
$
   K_w(z) \in A^p.
$ %\frac 1{(1-\langle z, w\rangle)^{n+1}} 
Now %for $g\in A^{p'}= \left( A^p\right)^*$, %with (\ref{dual-a}) 
%from Lemma   
we have 
$$ %\beqm
 \sup_{g\in B_{A^{p'}} }\left (\int_{\B_n} |\langle K_w, g \rangle |^{p'}  dv(w)  \right) ^{\frac 1{p'}} = \sup_{g\in B_{A^{p'}} }\left (\int_{\B_n} |g(w) |^{p'} dv(w)  \right) ^{\frac 1{p'}}=1.
$$% \eqm
We define a continuous mapping  $F$ on $\B_n $  as
\beqm 
 F: \B_n &\to & A^{p} \\
w &\mapsto & K_w(\cdot)
\eqm% 
Then, applying Lemma \ref{lem8}   to get
  \beqm
     \int_{\B_n} \|K(\cdot, w)\|_{L^2(\mu)}^{p'} dv(w) &=& \int_{\B_n} \|J_\mu K_w \|_{L^2(\mu)}^{p'} dv(w)\\
   &\le & \left(\pi_2 (J_\mu)\right)^{p'} \sup_{g\in B_{A^{p'}} } \int_{\B_n} |\langle K_w, g \rangle |^{p'}  dv(w)  .
  \eqm
Therefore,
\begin{equation} \label{p-small}
      \int_{\B_n} \|K(\cdot, w)\|_{L^2(\mu)}^{p'} dv(w) \le  \left(\pi_2 (J_\mu)\right)^{p'}.
\end{equation}

Conversely, suppose the statement (2) holds.  Define an operator $R$ on $L^2(\mu)$  as
$$
    R(g)(z)= \int_{\B_n} \frac {g(w)}{(1-\langle z, w\rangle)^{n+1}} d\mu(w).
$$
Then we have 
\beqm
    \|Rg\|_{L^{p'}}^{p'} &\le& \int_{\B_n}\left( \int_{\B_n} \frac {|g(w)|}{|1-\langle z, w\rangle|^{n+1}} d\mu(w) \right)^{p'}dv(z)\\
    &\le & \|g\|_{L^2(\mu)} ^{p'}   \int_{\B_n} \|K_z\|_{L^2(\mu)} ^{p'} dv(z)<\infty. 
\eqm
This means that $R$ is bounded from $L^2(\mu)$ to $L^{p'}$. 
Now for $f\in A^p$ and $g\in L^2(\mu)$,  
\beqm
  &&\int_{\B_n} |f(z)|\left ( \int_{\B_n}  \frac {|g(w)|}{|1-\langle z, w\rangle|^{n+1}} d\mu(w)) \right) dv(z)\\
  &\le & \|g\|_{L^2(\mu)}\int_{\B_n} |f(z)|\|K_z\|_{L^2(\mu)}  dv(z)\\
   &\le & \|g\|_{L^2(\mu)} \|f\|_{A^p} \left( \int_{\B_n} \|K_z\|_{L^2(\mu)}  ^{p'} dv(z)\right)^{\frac 1 {p'}}
  <\infty.
\eqm
This allows us to apply Fubini's theorem and then to use Bergman projection, we have
$$
    \langle J_\mu f , g\rangle_{L^2 (\mu)} = \langle  f , Rg\rangle_{A^2}.
$$
Therefore, $J_\mu: A^p\to L^2(\mu)$ is bounded and  $J_\mu^{*} =R$.
 Something more, it is trivial to verify that $R: L^2(\mu)\to L^{p'}$ is order bounded. Lemma \ref{lem66}  tells us that  $R$  is $p'$-summing. Since  $R(L^2(\mu))\subset A^{p'}$, we have   $R\in \Pi_{p'}(L^2(\mu), A^{p'})$ with
$$
          \pi_{p'} (R) \lesssim \left \| \sup_{g\in B_{L^2(\mu)}} |R(g)|\right\|_{A^{p'}} \lesssim \left ( \int_{\B_n} \|K(\cdot, w)\|_{L^2(\mu)}^{p'} dv(w) \right)^{\frac 1{p'}}.
$$
 This together with  \cite[Theorem 2.21]{DJT95} shows  $J_\mu: A^p \to L^2(\mu)$ is 2-summing and
\begin{equation}\label{R-b}
    \pi_2\left ({J_\mu }\right) \lesssim \left ( \int_{\B_n} \|K(\cdot, w)\|_{L^2(\mu)}^{p'} dv(w) \right)^{\frac 1{p'}}.
\end{equation}   
    
The norm equivalence (\ref{small-b}) follows from (\ref{p-small}) and (\ref{R-b}). 
\end{proof}  
\vspace{2mm}

\textit{Proof  of Theorem \ref{THM1}}. Recall that $ s=\frac{2p}{2p-2q+pq} $. 
First, we prove  that 
\begin{equation}\label{norm-z}
\left\|\widehat{\mu}_\delta(\cdot)K(\cdot, \cdot)^{\frac{q}{2}}\right\|_{L^s}^{\fr1q}
 \simeq 
 \left\|\widetilde{\mu} (\cdot)K(\cdot, \cdot)^{\frac{q}{2}}\right\|_{L^s}^{\fr1q}.
 \end{equation}
That the left hand side can be dominated by the right hand side is an easy consequence of the fact that $\widehat{\mu}_\delta(z)\le C \widetilde{\mu} (z)$ for all $z\in \B_n$. To get the other direction estimate, for any fixed $\delta>0 $ we have
\beqm
   & & \widetilde{\mu}(z)K(z, z)^{\frac q 2} \\
   &\lesssim& K(z, z)^{\frac q 2}\int _{\B_n} \frac {|K(z,w)|^2}{K(z, z)} \widehat{\mu}_\delta(w) dv (w)  \\
    &=& (1-|z|^2)^{(n+1)\left (1-\frac q 2\right)} \int _{\B_n} \frac {(1-|w|^2)^{ (n+1) \frac q 2}}{|1-\langle z, w \rangle |^{2(n+1)}} \left ( \widehat{\mu}_\delta(w)  (1-|w|^2)^{- (n+1)  \frac q 2 }  \right) dv (w)\\
    &=& S_{ (n+1)\left (1-\frac q 2\right), \  (n+1) \frac q 2, \  2(n+1)    } \left ( \widehat{\mu}_\delta(\cdot) K(\cdot, \cdot)^{\frac q 2}   \right) (z). 
\eqm
Notice that  $s\ge 1$ and $s=1$ if and only if $p=2$.  Now for  $p<2$ applying Theorem A  from \cite{ZZ22} (or for $p=2$ applying  Theorem B  from \cite{ZZ22}) to obtain   
$$
\left\|\widetilde{\mu} (\cdot)K(\cdot, \cdot)^{\frac{q}{2}}\right\|_{L^s}^{\fr1q}
 \lesssim
 \left\|\widehat{\mu}_\delta(\cdot)K(\cdot, \cdot)^{\frac{q}{2}}\right\|_{L^s}^{\fr1q}
.
$$
Then we have the estimates (\ref{norm-z}). 

 Now we prove that \(J_\mu: A^p \rightarrow  L^q(\mu)\) is \(r\)-summing if and only if \(\widehat{\mu}_\delta(\cdot)K(\cdot, \cdot)^{\frac{q}{2}} \in  {L}^s\), and 
\begin{equation}\label{norm-a}
\pi_r(J_\mu: A^p \rightarrow  A^p) \simeq  \left\|\widehat{\mu}_\delta(\cdot)K(\cdot, \cdot)^{\frac{q}{2}}\right\|_{L^s}^{\fr1q}
 \end{equation}
for some (or any) $\delta>0$.\\

\textbf{Case 1: $1< p \le2$, $1\le q \le2$.}
Suppose $\widehat{\mu}_\delta(\cdot)K(\cdot, \cdot)^{\frac{q}{2}}\in L^s$ for some $\delta>0$. Write
$$
\omega(z)=\widehat{\mu}_\delta(z)^{\frac{2s}{p'}}K(z, z)^{\frac{sq}{p'}-1}, \ \ z\in \B_n.
$$
Since
$$
\big(\omega(z)K(z, z)\big)^{\frac{p'}{2}}=\left(\widehat{\mu}_\delta(z)K(z, z)^{\frac{q}{2}}\right)^s,
$$
we get $\omega(\cdot)K(\cdot, \cdot)\in L^s$. It follows from Lemma \ref{small-a} that $J_\omega: A^p\to L^2(\omega dv)$ is 2-summing. Moreover,
\begin{equation}\label{norm31}
\pi_2(J_\omega)\simeq \left\|\widehat{\mu}_\delta(\cdot)K(\cdot, \cdot)^{\frac{q}{2}}\right\|_{L^s}^{s\cdot\fr1{p'}}.
\end{equation}
On the other hand, it is easy to check
$$
\widehat{\mu}_\delta(z)\omega(z)^{-\frac{q}{2}}=\widehat{\mu}_\delta(z)^{1-\frac{sq}{p'}}K(z, z)^{\frac{q}{2}-\frac{sq^2}{2p'}}.
$$
From this we get
\beqm
\int_{\B_n}\left(\widehat{\mu}_\delta(z)\omega(z)^{-\frac{q}{2}}\right)^{\frac{2}{2-q}}dv(z)&=& \int_{\B_n}\widehat{\mu}_\delta(z)^sK(z, z)^{\frac{sq}{2}}dv(z)\\
&=&\left\|\widehat{\mu}_\delta(\cdot)K(\cdot, \cdot)^{\frac{q}{2}}\right\|_{L^s}^{s}.
\eqm
Set $A^2(\omega dv)=L^2(\omega dv)\cap H(\B_n)$. We claim that
\beqm 
 {\rm Id}_1: A^2(\omega dv)&\to& L^q(\mu) \\
f &\mapsto & f
\eqm
is bounded. In fact, for $f\in A^2(\omega dv)$, set 
$$
E=\{z\in \B_n: \omega(z)=0\},
$$
then it is trivial to see $\omega(z)=0$ if and only if $\widehat{\mu}_\delta(z)=0$, and 
$$
\int_{\B_n }|f(z)|^q\widehat{\mu}_\delta(z)dv(z)=\int_{\B_n\setminus E }|f(z)|^q\widehat{\mu}_\delta(z)dv(z).
$$
Lemma 51 in \cite{ZZ08} and H\"{o}lder's inequality yield
\beqm
\int_{\B_n}|f(z)|^qd\mu(z)&\le& \int_{\B_n }|f(z)|^q\widehat{\mu}_\delta(z)dv(z)\\
&\le&\left(\int_{\B_n\setminus  E}|f(z)|^2\omega(z)dv(z)\right)^{\fr q2}\left[\int_{\B_n\setminus E}\left(\widehat{\mu}_\delta(z)\omega(z)^{-\frac{q}{2}}\right)^{\frac{2}{2-q}}dv(z)\right]^{\fr {2-q}2}\\
&=&\left\|f\right\|_{A^2(\omega dv)}^{q}\cdot\left\|\widehat{\mu}_\delta(\cdot)K(\cdot, \cdot)^{\frac{q}{2}}\right\|_{L^s}^{\fr {s(2-q)}2}.
\eqm
This shows ${\rm Id}_1: A^2(\omega dv)\to L^q(\mu)$ is bounded with
\begin{equation}\label{norm32}
\|{\rm Id}_1\|_{A^2(\omega dv)\to L^q(\mu)}\le \left\|\widehat{\mu}_\delta(\cdot)K(\cdot, \cdot)^{\frac{q}{2}}\right\|_{L^s}^{\fr {s(2-q)}{2q}}.
\end{equation}
Consider the chain of maps:
$$
A^p\stackrel{J_\omega}{\longrightarrow}  L^2(\omega dv)\stackrel{{\rm Id}_1}{\longrightarrow} L^q(\mu),
$$
The estimate (\ref{ideal-a}) indicates $J_\mu\left(A^p\to L^q(\mu)\right)={\rm Id}_1\circ J_\omega$ is 2-summing. It follows from (\ref{cotype}) that $J_\mu: A^p\to L^q(\mu)$ is $r$-summing.  Furthermore, by (\ref{norm31}) and (\ref{norm32}) we get
\begin{equation}\label{norm33}
\pi_r(J_\mu)\simeq \pi_2(J_\mu)\lesssim \left\|\widehat{\mu}_\delta(\cdot)K(\cdot, \cdot)^{\frac{q}{2}}\right\|_{L^s}^{\fr1q}.
\end{equation}

Conversely, suppose $J_\mu\in \Pi_r( A^p, L^q(\mu))$. Since \(A^{p}\) and $L^q(\mu)$ both have cotype 2 when $1< p \le2$, $1\le q \le2$,  we have $\Pi_r( A^p, L^q(\mu))=\Pi_q( A^p, L^q(\mu))$. Pietsch domination theorem tells us there exists some regular probability measure $\sigma$ on ${B}_{A^{p'}}$ such that for each $f\in A^{p}$, there holds
 \begin{equation}\label{ineq-1}
\left\|f\right\|_{L^q(\mu)}^q\leq  \big({\pi }_{q}(J_\mu)\big)^q\cdot\int_{  B_{A^{p'}}}{\left|\left\langle g, f \right\rangle\right| }^q{d\sigma(g)}.
 \end{equation}
Now for $\{c_j\}_{j=1}^\infty$, set
 \begin{equation}\label{ft-def}
 f_t(z)=\sum_{j=1}^mc_jr_j(t)k_{a_j, p}(z),
  \end{equation}
 where  $r_j(t)$ is the Rademacher function on $[0,1]$, $a_j$ belongs to the $\delta$-lattice $\{a_j\}$ for small enough $\delta$. On one hand, since for $w\in B(z, \delta)$,
 $$
 |1-\langle z, w \rangle|\simeq 1-|z|^2\simeq 1-|w|^2, v\left(B(z, \delta)\right)\simeq K(z,z)^{-1},
 $$
Fubini's theorem and Khinchine's inequality give
\beqm
\int_{0}^{1}\left\|f_t\right\|_{L^q(\mu)}^{q}dt
&=&\int_{\B_n}\left(\int_{0}^{1}\left|\sum^{m}_{j=1}c_jr_j(t)k_{a_j, p}(z)\right|^{q}dt\right)d\mu(z)\\
&\simeq&\int_{\B_n}\left(\sum^{m}_{j=1}\left|c_j\right|^2|k_{a_j, p}(z)|^2\right)^{\fr q 2}d\mu(z)\\
&\gtrsim&\sum^{m}_{k=1}\left|c_k\right|^q\int_{B(a_k, \delta)}|k_{a_k, p}(z)|^qd\mu(z)\\
&\gtrsim&\sum^{m}_{k=1}|c_k|^q\fr{K(a_k, a_k)^q}{\|K(\cdot, a_k)\|^q_{A^p}}\cdot\mu(B(a_k, \delta))\\
&\simeq&\sum^{m}_{k=1}\left|c_k\right|^{q}K(a_k, a_k)^{\fr qp-1}\cdot\widehat{\mu}_\delta(a_k).
\eqm
On the other hand, thanks to the reproducing kernel formula, by using H\"{o}lder's inequality with  $\fr {p'} 2$ and its conjugate exponent, we have
\beqm
&&\int_{0}^{1}\int_{B_{A^{p'}}}\left|\left\langle g, f_t \right\rangle \right|^{q}{d\sigma(g)}dt\\
&=&\int_{B_{A^{p'}}}d\sigma(g)\int_{0}^{1}\left|\sum_{j=1}^mc_jr_j(t)\frac{g(a_j)}{\|K(\cdot, z)\|_{A^{p}}}\right|^{q}dt\\
&\simeq&\int_{B_{A^{p'}}}\left(\sum_{j=1}^m|c_j|^2|g(a_j)|^2K(a_j, a_j)^{-\fr2{p'}}\right)^{\fr q 2}d\sigma(g)\\
&\lesssim& \int_{B_{A^{p'}}}\left(\sum_{j=1}^m|c_j|^{\frac{2p'}{p'-2}}\right)^{\frac{p'-2}{p'}\cdot \fr q 2}\cdot\left(\sum_{j=1}^m |g(a_j)|^{p'}K(a_j, a_j)^{-1}\right)^{\frac{q}{2}\cdot\frac{2}{p'}}d\sigma(g)\\
&\lesssim& \|\{|c_j|^q\}\|_{l^{\frac{2p}{(2-p)q}}}\cdot\int_{B_{A^{p'}}}\|g\|_{A^{p'}}^qd\sigma(g)\\
&\le& \|\{|c_j|^q\}\|_{l^{\frac{2p}{(2-p)q}}}.
\eqm
Therefore, (\ref{ineq-1}) yields
$$
\sum^{m}_{j=1}\left|c_j\right|^{q}K(a_j, a_j)^{\fr qp-1}\widehat{\mu}_\delta(a_j)\lesssim \big({\pi }_{q}(J_\mu)\big)^q\cdot \|\{|c_j|^q\}\|_{l^{\frac{2p}{(2-p)q}}}.
$$
Since $m$ is arbitrary, we get
$$
\sum^{\infty}_{j=1}\left|c_j\right|^{q}K(z_j, z_j)^{\fr qp-1}\widehat{\mu}_\delta(z_j)\lesssim \big({\pi }_{q}(J_\mu)\big)^q\cdot \|\{|c_j|^q\}\|_{l^{\frac{2p}{(2-p)q}}}.
$$
A duality argument to $\frac{2p}{(2-p)q}$ and its conjugate $s$, notice $(\fr qp-1)\cdot s=\fr q2\cdot s-1$, we obtain
$$
\sum^{\infty}_{j=1}\left(\widehat{\mu}_\delta(z_j)K(z_j, z_j)^{\fr q2}\right)^s\cdot |B(z_j, \delta)|\lesssim \big({\pi }_{q}(J_\mu)\big)^{qs}.$$
This indicates $\widehat{\mu}_\delta(\cdot)K(\cdot, \cdot)^{\frac{q}{2}}\in L^s$. Moreover,
\begin{equation}\label{norm36}
\left\|\widehat{\mu}_\delta(\cdot)K(\cdot, \cdot)^{\frac{q}{2}}\right\|_{L^s}^{\fr1q}\lesssim \pi_q(J_\mu)\simeq \pi_r(J_\mu).
\end{equation}\\

\textbf{Case 2: $p=1$, $1\le q \le2$.}
Suppose $\widehat{\mu}_\delta(\cdot)K(\cdot, \cdot)^{\frac{q}{2}}\in L^{\fr 2{2-q}}$, we show ${\rm Id}_2: A^1\to A^2_2$ is 1-summing, and ${\rm Id}_3: A^2_2\to L^q(\mu)$ is bounded.
In fact, for $f\in A^1$, it follows from (\ref{formula1}) we get
\beqm
\|f\|_{A_2^2}^2&=&\int_{\B_n}(1-|z|^2)^2|f(z)|^2dv(z)\\
&\le& \|f\|_{A^1}\cdot\int_{\B_n}|f(z)|dv(z)\\
&=&\|f\|_{A^1}^2,
\eqm
which gives ${\rm Id}_2: A^1\to A^2_2$ is bounded with $\|{\rm Id}_2\|_{A^1\to A^2_2}\le 1$. It follows from \cite[Theorem 3.4]{DJT95} that ${\rm Id}_2: A^1\to A^2_2$ is 1-summing with
 \begin{equation}\label{norm-3-7}
\pi_1({\rm Id}_2: A^1\to A^2_2)\lesssim \|{\rm Id}_2\|_{A^1\to A^2_2}\le 1.
\end{equation}
Meanwhile, for $f\in A^2_2$, we get
\beqm
\|{\rm Id}_3f\|^q_{L^q(\mu)}&=&\|f\|^q_{L^q(\mu)}\\
&\lesssim& \int_{\B_n}|f(\xi)|^q\cdot\widehat{\mu}_\delta(\xi)dv(\xi)\\
&=& \int_{\B_n}\big[(1-|\xi|^2)^q|f(\xi)|^q\big]\cdot\left[\widehat{\mu}_\delta(\xi)K(\xi,\xi)^{\frac{q}{2}}\right]dv(\xi)\\
&\le & \left(\int_{\B_n}(1-|\xi|^2)^2|f(\xi)|^2dv(\xi)\right)^{\fr q2}\cdot\left(\int_{\B_n}\widehat{\mu}_\delta(\xi)^{\frac{2}{2-q}}\cdot K(\xi,\xi)^{\frac{q}{2-q}} dv(\xi)\right)^{\fr{2-q}2}\\
&=& \|f\|_{A_2^2}^q\cdot \left\|\widehat{\mu}_\delta(\cdot) K(\cdot,\cdot)^{\frac{q}{2}}\right\|_{L^ {\fr 2{2-q}}}.
\eqm
This gives ${\rm Id}_3: A^2_2\to L^q(\mu)$ is bounded with
 \begin{equation}\label{norm-3-8}
\|{\rm Id}_3\|_{A^2_2\to L^q(\mu)}\lesssim \left\|\widehat{\mu}_\delta(\cdot) K(\cdot,\cdot)^{\frac{q}{2}}\right\|^{\fr 1q}_{L^ {\fr 2{2-q}}}.
\end{equation}
Consider the chain of maps:
$$
A^1\stackrel{{\rm Id}_2}{\longrightarrow}  A^2_2\stackrel{{\rm Id}_3}{\longrightarrow} L^q(\mu),
$$
we get $J_\mu(A^1\to L^q(\mu))={\rm Id}_3\circ {\rm Id}_2$ is 1-summing. Furthermore, by (\ref{norm-3-7})and (\ref{norm-3-8}) we get
\begin{equation}\label{norm39}
\pi_1\big(J_\mu: A^1\to L^q(\mu)\big)\lesssim \left\|\widehat{\mu}_\delta(\cdot)K(\cdot, \cdot)^{\frac{q}{2}}\right\|_{L^{\fr 2{2-q}}}^{\fr1q}.
\end{equation}

Conversely, suppose $J_\mu: A^1\to L^q(\mu)$ is 1-summing, then $J_\mu\in \Pi_q(A^1, L^q(\mu))$. By Pietsch domination theorem, there exists regular probability measure $\sigma$ on ${B}_{\mathcal{B}}$ such that for each $f\in A^{1}$, there holds
 \begin{equation}\label{3-10}
\left\|f\right\|^q_{L^{q}(\mu)}\leq \left ({\pi }_{q}(J_\mu)\right)^q\cdot\int_{  B_{\mathcal{B}}}{\left|\left\langle g, f \right\rangle\right|^q }{d\sigma(g)}.
 \end{equation}
For $z, w \in \B_n$, set
$$
L(z, w)=\frac{\partial}{\partial z_1}K(z, w)=\frac{(n+1)\overline{w}_1}{(1-\langle z, w\rangle)^{n+2}}.
$$
Then  for $g\in \mathcal{B}\subset A^2$,
by $g(z)= \langle g, K(\cdot, z) \rangle $
we know
$$
\frac {\partial}{\partial z_1} g\left (z_1, z_2, \cdots, z_n\right )%=\frac{d}{d z}\langle g, K(\cdot, z) \rangle
=\left\langle g, \overline{\frac{\partial}{\partial z_1 }K(z ,  \cdot)} \right\rangle= \langle g, \overline{L(z, \cdot) }\rangle .
$$
Meanwhile, applying   \cite[Theorem 1.12]{Zh05}
to obtain
$$
\|L(\cdot, w)\|_{A^1}\simeq \int_{\B_n}\frac{1}{|1-\langle w, z\rangle|^{n+2}}dv(z)\simeq \frac{1}{1-|w|^2},
$$
and
$$
\int_{B(w, \delta)}\frac{|L(z, w)|^q}{\|L(\cdot, w)\|^q_{A^1}}d\mu(z)\simeq \int_{B(w, \delta)}\frac{1}{(1-|w|^2)^{
(n+1)q}}d\mu(z)\simeq K(w, w)^{q-1}\widehat{\mu}_\delta(w).
$$
Now for $\{c_j\}_{j=1}^\infty$, take
 $$f_t(z)=\sum_{j=1}^mc_jr_j(t)\frac{L(z, a_j)}{\|L(\cdot, a_j)\|_{A^1}},$$
 where  $r_j(t)$ is the Rademacher function on $[0,1]$,   $\{a_j\}$ is any given $\delta$-lattice with  $\delta>0$.
Fubini's theorem and Khinchine's inequality give
\beqm
\int_{0}^{1}\left\|f_t\right\|^q_{L^{q}(\mu)}dt
&\simeq&\int_{\B_n}\left(\sum^{m}_{k=1}\left|c_k\right|^2\frac{|L(z, a_k)|^2}{\|L(\cdot, a_k)\|^2_{A^1}}\right)^{\fr q 2}d\mu(z)\\
&\gtrsim&\sum^{m}_{j=1}\left|c_j\right|^q\int_{B(a_j, \delta)}\frac{|L(z, a_j)|^q}{\|L(\cdot, a_j)\|^q_{A^1}}d\mu(z)\\
&\gtrsim& \sum_{j=1}^m |c_j|^qK(a_j, a_j)^{q-1}\widehat{\mu}_\delta(a_j).
\eqm
On the other hand, by $$\sup_{j}  \left |\frac{\partial g}{\partial z_1} (a_j)\right |(1-|a_j|^2)\le \|g\|_{\mathcal{B}},$$ we get
\beqm
\int_{0}^{1}\int_{B_{\mathcal{B}}}{\left|\left\langle g, f_t \right\rangle\right|^q }{d\sigma(g)}dt
&\simeq&\int_{B_{\mathcal{B}}}\left(\sum_{j=1}^m|c_j|^2\left |\frac{\partial g}{\partial z_1} (a_j)\right |^2(1-|a_j|^2)^2\right)^{\fr q 2}d\sigma(g)\\
&\lesssim& \left(\sum_{j=1}^m|c_j|^{2}\right)^{\frac{q}{2}}\cdot \int_{B_{\mathcal{B}}}\|g\|^{q}_{\mathcal{B}}d\sigma(g)\\
&\le& \|\{|c_j|^q\}\|_{l^{\frac{2}{q}}}.
\eqm
Therefore, (\ref{3-10}) deduces that
$$
\sum^{\infty}_{j=1}|c_j|^qK(a_j, a_j)^{q-1}\widehat{\mu}_\delta(a_j)\le\left( {\pi }_{q}(J_\mu)\right)^q\cdot \|\{|c_j|^q\}\|_{l^{\frac{2}{q}}}.
$$
 A duality argument to the exponentials  $\frac{2}{q}$ and its conjugate $\frac{2}{2-q}$, we obtain
$$
\sum^{\infty}_{j=1}\left[\widehat{\mu}_\delta(z_j)K(a_j, a_j)^{\fr q {2}}\right]^{\fr 2 {2-q}}\cdot |B(a_j, \delta)|\lesssim \big({\pi }_{q}(J_\mu)\big)^{\frac{2q}{2-q}}.
$$
This shows $\widehat{\mu}_\delta(\cdot)K(\cdot, \cdot)^{\frac{q}{2}}\in L^{\frac{2}{2-q}}$ with the estimate
\begin{equation}\label{norm311}
\left\|\widehat{\mu}_\delta(\cdot)K(\cdot, \cdot)^{\frac{q}{2}}\right\|_{L^{\frac{2}{2-q}}}^{\fr1q}\lesssim \pi_q\big(J_\mu: A^1\to L^q(\mu)\big)\simeq \pi_r\big(J_\mu: A^1\to L^q(\mu)\big).
\end{equation}

Combine the proof above we know that for $1\le p, q\le 2$  and  \(r \geq  1\), \(J_\mu: A^p \rightarrow  L^q(\mu)\) is \(r\)-summing if and only if \(\widehat{\mu}_\delta(\cdot)K(\cdot, \cdot)^{\frac{q}{2}} \in  {L}^s\) for some (of any) $\delta>0$.
And the desired norm equivalence (\ref{norm121})  comes from (\ref{norm33}), (\ref{norm36}), (\ref{norm39}) and (\ref{norm311}).
The proof of Theorem \ref{THM1} is complete.
%\end{proof}

\section{Proof of Theorem \ref{THM2}}

In this section, we are going to present the proof of Theorem \ref{THM2}. Our proof will be carried out  in four subsections.

\subsection{The $r$-summing Toeplitz operators on $A^{p}$, $1< p \le2$.}

\begin{proposition}\label{prop3.1}
 Suppose $1< p \le  2$, $1< q \le  2$ and $r\ge 1$, $s=\frac{2p}{2p-2q+pq}$. Let $\mu$ be a positive Borel measure on \(\B_n\). Then  \({T}_{\mu} : A^{p} \rightarrow  A^{q}\) is \(r\)-summing if and only if  \(\widehat{\mu}_\delta(\cdot)K(\cdot, \cdot)^{\fr 12} \in L^ {qs}\). Furthermore,
\begin{equation}\label{norm41}
\pi_r(T_\mu: A^p\to A^q)\simeq\left\|\widehat{\mu}_\delta(\cdot)K(\cdot, \cdot)^{\frac{1}{2}}\right\|_{L^{qs}}.
 \end{equation}
\end{proposition}

\begin{proof} Suppose $T_\mu\in \Pi_r(A^{p}, A^q)$, \(1< p \le  2, 1\ge q \le  2\) and $1\le r<\infty$.  Since both \(A^{p}\) and \(A^{q}\)  have cotype 2,  we know  $ \Pi_r(A^{p}, A^q)=\Pi_q(A^{p}, A^q)$.   Pietsch domination theorem tells us there exists regular probability measure $\sigma$ on ${B}_{A^{p'}}$ such that, for  $f\in A^{p}$, there holds
 \begin{equation}\label{42}
\left\|T_\mu\left(f\right)\right\|_{A^{q}}^{q}\leq  \big({\pi }_{q}\left(T_\mu \right)\big)^{q}\cdot\int_{  B_{A^{p'}}}{\left|\left\langle g, f \right\rangle\right| }^{q}{d\sigma(g)}.
 \end{equation}
Some element calculate shows $ T_\mu(k_{z, p})(z)=\widetilde{\mu}(z)\cdot K(z, z)^{\fr1p}$. Hence,
\beqm
\widetilde{\mu}(z)^q&=&|T_\mu(k_{z, p})(z)|^qK(z, z)^{-1}K(z, z)^{1-\fr qp}\\
&\lesssim& K(z, z)^{1-\fr qp}\int_{B(z, \delta)}|T_\mu(k_{z, p})(z)|^qdv(z).
\eqm
Now fix some $\delta$-lattice $\{a_j\}$. For $\{c_j\}_{j=1}^\infty$, take $f_t$ as in (\ref{ft-def}).
Fubini's theorem and Khinchine's inequality give
\beqm
\int_{0}^{1}\left\|T_\mu\left(f_t\right)\right\|_{A^{q}}^{q}dt
&\simeq&\int_{\B_n}\left(\sum^{m}_{k=1}\left|c_k\right|^2|T_\mu(k_{a_k, p})(z)|^2\right)^{\fr q 2}dv(z)\\
&\gtrsim&\sum^{m}_{j=1}\left|c_j\right|^q\int_{B(a_j, \delta)}|T_\mu(k_{a_j, p})(z)|^qdv(z)\\
&\gtrsim&\sum^{m}_{j=1}\left|c_j\right|^{q} K(a_j, a_j)^{\fr qp-1}\widetilde{\mu}(a_j)^q
\eqm
On the other hand, as in Section 3.1, we get
\beqm
\int_{0}^{1}\int_{B_{A^{p'}}}\left|\left\langle g, f_t \right\rangle \right|^{q}{d\sigma(g)}dt
\lesssim  \|\{|c_j|^q\}\|_{l^{\frac{2p}{(2-p)q}}}.
\eqm
With (\ref{42}) we obtain
$$
\sum^{\infty}_{j=1}\left|c_j\right|^{q} K(a_j, a_j)^{\fr qp-1}\widetilde{\mu}(a_j)^q\lesssim \left(\pi_q( T_\mu) \right)^{q}\cdot\|\{|c_j|^q\}\|_{l^{\frac{2p}{(2-p)q}}}.
$$
A duality argument  shows
$$
\sum^{\infty}_{j=1}\left(\widetilde{\mu}(z_j)^q K(a_j, a_j)^{\fr q2}\right)^sK(a_j, a_j)^{-1}\lesssim \left(\pi_q( T_\mu) \right)^{qs}.
$$
This indicates $\widehat{\mu}_\delta(\cdot)K(\cdot, \cdot)^{\frac{1}{2}}\in L^{qs}$ with the following estimate
\begin{equation}\label{norm43}
\left\|\widehat{\mu}_\delta(\cdot)K(\cdot, \cdot)^{\frac{1}{2}}\right\|_{L^{qs}}\lesssim \pi_q\big(T_{\mu}: A^p\to A^q\big)\simeq \pi_r\big(T_{\mu}: A^p\to A^q\big).
\end{equation}

Conversely, suppose $\widehat{\mu}_\delta(\cdot)K(\cdot, \cdot)^{\fr 12} \in L^ {qs}$,  \(1<p \le  2, 1<q\le 2\).  It is enough to show $T_\mu\in \Pi_2(A^{p}, A^q)$ .
For this purpose we set  $\omega(z)=\widehat{\mu}_\delta(z)^q$.
For $w\in B(z, \delta/2)$, by  $B(z, \delta/2)\subset B(w, \delta)\subset B(z, 2\delta)$  we know
$$
\widehat{\mu}_\frac{\delta}{2}(z)\lesssim \widehat{\mu}_\delta(w)\lesssim \widehat{\mu}_{2\delta}(z).
$$
%taking the $q$-th power, integrating over $D(z,\delta)$, and then dividing by $|D(z,\delta)|$, we get
It implies
$$
\widehat{\mu}_\frac{\delta}{2}(z)^q\lesssim \widehat{\omega}_\delta(z)\lesssim \widehat{\mu}_{2\delta}(z)^q.
$$
Notice that $\omega(\cdot)K(\cdot, \cdot)^{\fr q2}=\left[\widehat{\mu}_\delta(\cdot)K(\cdot, \cdot)^{\fr 12}\right]^q$  and $\widehat{\mu}_\delta(\cdot)K(\cdot, \cdot)^{\fr 12} \in L^ {qs}$ is independent of the choice of $\delta$, we know $\widehat{\omega}_\delta(\cdot)K(\cdot, \cdot)^{\fr q2} \in L^ {s}$.  Theorem \ref{THM1} tells us
that $J_\omega: A^p\to L^q(\omega dv),$ $f\mapsto f$ is $r$-summing. Set $A^q(\omega dv)=L^q(\omega dv)\cap H(\B_n)$, it is easy to check that $A^q(\omega dv)$ is a Banach space. As an operator from $A^p$ to $A^q(\omega dv),$ $J_\omega$ is also bounded. Moreover,
\begin{equation}\label{L-1}
\pi_r\big(J_\omega:  A^p\to A^q(\omega dv)\big)\simeq\left\|\widehat{\mu}_\delta(\cdot)K(\cdot, \cdot)^{\frac{1}{2}}\right\|_{L^{qs}}.
\end{equation}
Define the operator $T$ from $A^q(\omega dv)$ to $L^q$ by
$$
T:f\mapsto \int_{\B_n} f(w)K(z, w)d\mu(w).
$$
We claim that $T : A^q(\omega dv)\to L^q$ is bounded. To see this, take the operator $S_{b,c}$  as
\begin{equation}\label{S_bc}
S_{b,c}f(z)=\int_{\B_n}f(w)\frac{(1-|w|^2)^{b}}{|1-\langle z, w\rangle|^{c}}dw(w).
\end{equation}
It is easy to see
$$
|T f(z)|\lesssim \int_{\B_n} |f(w)||K(z, w)|\widehat{\mu}_\delta(w)dv(w)= S_{0, n+1} (|f|\widehat{\mu}_\delta)(z).
$$
Since $1< q<\infty$, it follows from \cite[Theorem 3]{Zhao15} that $S_{ 0, n+1}$ is bounded on $L^q$.
 Therefore, for $f\in A^q(\omega dv)$,
\begin{equation}\label{L-2}
 \|T f\|_{L^q}\lesssim  \|S_{0, n+1}\|_{L^q\to L^q}\|f\cdot\widehat{\mu}_\delta\|_{L^q}=\|S_{0, n+1}\|_{L^q\to L^q}\|f\|_{A^q(\omega dv)}.
\end{equation}
This indicates $T: A^q(\omega dv)\to L^q$ is bounded with the norm estimate $\|T\|_{A^q(\omega dv)\to L^q}\lesssim 1$.
Now let us see
the following commutative diagram:
$$
			\xymatrix{
				\  \ \ \ \ A^{p} \ \ \ \ 	\ar[r]^{\ T_\mu}\ \  \ar[d]^{J_\omega} & \ \ \ \ A^{q} \ \ \ \ \\
				 \ \ \ \ A^{q}(\omega dv) \ \ \ \ 	\ar[r]_{\ \ \ T} & \ \ \ \ L^{q} \ \ \ \ , \ar[u]_{P}}
 $$
where $P$ is the usual Bergman projection as (\ref {projection}). 
%Choose the parameter $\gamma$ such that $p(\gamma+1)>1$, \cite[Theorem 4.24]{Zh07}  \textcolor[rgb]{0.00,0.07,1.00}{\cite{Zh05}??} shows that for
For $1< p<\infty$, $P$ is bounded from $L^p$ to $A^p$ with $\|P\|_{ L^p\to A^p}\lesssim 1$ and $P(f)=f$ for $f\in A^p$. Therefore, (\ref{ideal-a}) and  (\ref{L-1}), (\ref{L-2}) turn out  $T_\mu=P \circ T\circ J_\omega$ is $r$-summing with
\beqm
\pi_r(T_\mu: A^{p}\to A^{q})&\le& \pi_r\left(J_\omega: A^p\to A^q(\omega dv)\right)\cdot\|T\|_{A^q(\omega dv)\to L^q}\cdot\|P\|_{L^p \rightarrow A^p}\\
&\lesssim& \left\|\widehat{\mu}_\delta(\cdot)K(\cdot, \cdot)^{\frac{1}{2}}\right\|_{L^{qs}}.
\eqm
This and the estimate  (\ref{norm43}) turn out  the norm equivalence (\ref{norm41}).
\end{proof}

Take $p=q$ in Proposition \ref{prop3.1}, we have the following conclusion  immediately.

\begin{corollary}\label{cor3.1}
 Suppose \(1< p \le  2\) and $r\ge 1$. Let $\mu$ be a positive Borel measure on \(\B_n\). Then  \({T}_{\mu} : A^{p} \rightarrow  A^{p}\) is \(r\)-summing if and only if  \(\widetilde{\mu} \in L^2(d\lambda)\). Furthermore,
$$
\pi_r(T_\mu)\simeq\left\|\widetilde{\mu}\right\|_{L^2(d\lambda)}.
$$
\end{corollary}

\subsection{The $r$-summing Toeplitz operators on $A^{p}$, $p\ge 2$ and $1\le r\le p'$.}\label{4.2}

% The purpose of this subsection is to prove the following proposition.

\begin{proposition}\label{prop3.5}
 Suppose $p\ge 2$ and $1\le r\le p'$. Let $\mu$ be a positive Borel measure on \(\B_n\). Then  \({T}_{\mu} : A^{p} \rightarrow  A^{p}\) is \(r\)-summing if and only if  \(\widetilde{\mu} \in L^ {p'}(d\lambda)\). Furthermore,
\begin{equation}\label{norm4-7}
\pi_r(T_\mu)\simeq\|\widetilde{\mu}\|_{L^{p'}(d\lambda)}.
 \end{equation}
\end{proposition}

\begin{proof}
  Suppose $T_\mu\in \Pi_{r}(A^{p}, A^p)$, since $p'\ge r$, it is clear that $T_\mu\in \Pi_{p'}(A^{p}, A^p)$ with $\pi_{p'}(T_\mu: A^{p} \rightarrow A^p)\le\pi_{r}(T_\mu: A^{p}\rightarrow A^p)$.  Meanwhile, since $\Pi_{p'}(A^{p}, A^p)\subset  B(A^{p}, A^p)$, we have $(T_\mu)^*= T_\mu \in B(A^{p'}, A^{p'})$, here $B(X, Y)$ is the set of bounded linear operators from $X$ to $Y$.  Consider
$$
L^p\stackrel{P}{\longrightarrow} A^p\stackrel{T_\mu}{\longrightarrow} A^p,\ \textrm{and} \ \ A^{p'}\stackrel{T_\mu}{\longrightarrow} A^{p'} \stackrel{\mathrm{Id}_1} {\longrightarrow} L^{p'},
$$
 Then $P$ is bounded from $L^p$ to $A^p$.
Set $Q_1 = T_\mu\circ P: A^{p}\to A^{p}$ and $Q_2 = {\rm Id}_1\circ {T}_{\mu}: A^{p'} \rightarrow  L^{p'}$, then $Q_1= Q_2^*$ and  $Q_1\in {\Pi_{p'}} (L^{p}, A^{p})$.
 It follows from \cite[Corollary 5.21]{DJT95} that \(Q_2: A^{p'} \rightarrow  L^{p'}\) is order bounded, that is $\mathop{\sup }\limits_{{g \in  {B}_{A^{p'}}}}{\left| T_\mu(g)(z) \right| }\in L^{p'}$ with the estimate
 \begin{equation}\label{46}
\big\|\mathop{\sup }\limits_{{g \in  {B}_{A^{p'}}}}{\left| T_\mu(g)\right| }\big\|_{L^{p'}}\le \pi_{p'}(T_\mu: A^{p}\to  A^p).
 \end{equation}
Notice that
$$
k_{z, p'}(\xi)=\frac{K(\xi, z)}{\|K(\cdot, z)\|_{A^{p'}}}\simeq \frac{K(\xi, z)}{K(z, z)^{1-\frac{1}{p'}}},
$$
we have
\beqm
\mathop{\sup }\limits_{{g \in  {B}_{A^{p'}}}}{\left| T_\mu(g)(z) \right| }
\ge |T_\mu(k_{z, p'})(z)|
\simeq \widetilde{\mu }(z)\cdot K(z, z)^{\frac{1}{p'}}.
\eqm
Hence, (\ref{46}) yields
$$
\int_{\B_n}\widetilde{\mu }(z)^{p'}K(z, z)dv(z)\lesssim\big(\pi_{p'}(T_\mu: A^{p}\to  A^p)\big)^{p'},
$$
which shows \(\widetilde{\mu} \in  {L}^{p'}(d\lambda)\) with
 \begin{equation}\label{4.3333}
 \|\widetilde{\mu}\|_{L^{p'}(d\lambda)}\lesssim\pi_{p'}(T_\mu: A^{p}\to  A^p)\le \pi_{r}(T_\mu: A^{p}\to  A^p).
 \end{equation}

Conversely, suppose \(\widetilde{\mu} \in L^ {p'}(d\lambda)\).
We  fix some  $N>0 $ so that
$$
\alpha= N+(n+1)\left( \frac 1p-1\right)>0
$$
and
$$
\beta= 2N+ (n+1) \left (1-\frac 2 p \right)> 0.
$$
We consider the chain of maps:
$$
A^p\stackrel{R^{0, N} \circ T_\mu}{\longrightarrow}A_{\alpha }^1\stackrel{{\rm Id}_2}{\longrightarrow}A^2_{\beta }\stackrel{{\rm Id}_3}{\longrightarrow}A^p_{Np}\stackrel{R_{0,N} }{\longrightarrow}A^p.
$$
We claim the following statements are true:

%(i)  $R^{0, N} \circ T_\mu: A^p\to A^p$ is 1-summing.

(i)  $R^{0, N}\circ T_\mu: A^p\to A_{\alpha}^1$ is bounded.

(ii) All of  ${\rm Id}_2$, ${\rm Id}_3$ and  $R_{0,N}$ are bounded.

(iii)  ${\rm Id}_2: A_{\alpha }^1\to A^2_{\beta }$ is 1-summing.

To prove (i), from Lemma \ref{R-op}  we see that
\beqm
   R^{0, N} T_\mu f(z)&=& R^{0, N} \int_{\B_n} f(w) K(z, w) d\mu(w) \\
   &=& \int_{\B_n} f(w) R^{0, N} \left (K(\cdot, w)\right)(z) d\mu(w)\\
   &=& \int_{\B_n} f(w) \frac 1{ (1-\langle z, w\rangle) ^{n+1+N}} d\mu(w).
\eqm
Then by \cite[Theorem 1.12]{Zh05}  we have
\beqm
     \left\| R^{0, N} T_\mu f\right\|_{A^1_{\alpha}} &\le &
     \int _{\B_n} (1-|z|^2)^\alpha dv(z) \int_{\B_n} \frac {|f(w)|\widehat{\mu}_\delta(w)}{|1-\langle z, w\rangle|^{n+1+N}}dv(w)\\
     &=& \int_{\B_n} |f(w)|\widehat{\mu}_\delta(w) dv(w) \int _{\B_n} \frac {(1-|z|^2)^\alpha}{|1-\langle z, w\rangle|^{n+1+N}}dv(z)\\
     &\lesssim&  \int_{\B_n} |f(w)|\widehat{\mu}_\delta(w) K(w, w)^{\frac 1{p'} } dv(w).
\eqm
Applying Fubini's Theorem and the Holder's inequality with the exponentials $p$ and $p'$ to get
 \begin{equation} \label{1<r<p'}
   \left\| R^{0, N} T_\mu f\right\|_{A^1_{\alpha}}  \lesssim \|f\|_{A^p} \|\widehat{\mu}_\delta \|_{L^{p'}(d\lambda)}.
 \end{equation}
This gives the boundedness of $R^{0, N}\circ T_\mu$ from $A^p$  to $A_{\alpha}^1$.

To see (ii), we know   from \cite[Exercise 2.27]{Zh05} that,  for $0<p<q<\infty$ and
$$
n+1+\alpha =\frac{n+1+\beta}{2},
$$
the inclusion  ${\rm Id}_2: A_{\alpha}^1\to A^2_{\beta}$ %and ${\rm Id}_3: A^2_{\fr2p(kp+2)-2}\to A^p_{kp}$ are
is bounded. Moreover,
\begin{equation} \label{fp-a}
\|{\rm Id}_2\|_{A_{\alpha }^1\to A^2_{\beta }}\lesssim 1. %\ \ \|{\rm Id}_2\|_{A^3_{\fr2p(kp+2)-2}\to A^p_{kp}}\lesssim 1.
\end{equation}
That ${\rm Id}_3$ is bounded with the estimate
\begin{equation} \label{fp-b}
%\|{\rm Id}_3\|_{A_{k+\fr2p-2}^1\to A^2_{\fr2p(kp+2)-2}}\lesssim 1.
\|{\rm Id}_3\|_{A^2_{\beta }\to A^p_{Np}}\lesssim 1
\end{equation}
can be seen easily from estimate (\ref{weighted-a}) on $A^2_\beta$.

Now  we show that  $R_{0,N} : A^p_{Np}\to A^p$ is bounded. For $f\in H^\infty,$
from \cite[Theorem 2.2]{Zh05} we know
$$
    f(z)= c_{N} \int_{\B_n} \frac {f(w)(1-|w|^2)^{N}}{(1-\langle z, w\rangle)^{n+1+N} } dv(w).
$$
Therefore,
\beqm
   R_{0,N} f(z) &=& \int_{\B_n} {f(w)(1-|w|^2)^{N }}R_{0,N} \left(\frac 1 {(1-\langle z, w\rangle)^{n+1+N } }\right)(w) dv(w)\\
   %&=&  \int_{\B} {f(w)} \frac {(1-|w|)^{N }} {(1-z\bar{w})^{n+1} } dv(w)\\
   &=&  \int_{\B_n} \left(f(w)(1-|w|^2)^N\right) \frac {1} {(1-\langle z, w\rangle)^{n+1} } dv(w)\\
    &=&  P\left((1-|\cdot|^2)^N f \right)(z).
\eqm
With the boundedness of the Bergman projection $P$ on $L^p$ and the density of $H^\infty$ in $A^p_{pN}$, we have
\begin{equation}\label{fp-c}
  \| R_{0, N} f\|_{A^p} \le \left \|P\left ((1-|\cdot|^2)^N f \right)\right\|_{A^p} \le C \|P\|_{L^p\to A^p} \|f \|_{A^p_{Np}}.
\end{equation}

Now we prove (iii).
% Set
%$$
%\textcolor[rgb]{0.00,0.07,1.00}{\alpha=k+\fr2p-2\ \ \textrm {  and} \ \  \beta=\fr2p(kp+2)-2}
%.$$
Consider the maps:
$$
L^1_\alpha\stackrel{P_\gamma}{\longrightarrow} A^1_\alpha\stackrel{{\rm Id}_2}{\longrightarrow} A^2_\beta\stackrel{{\rm Id}_4}{\longrightarrow}L^2_\beta,
$$
where 
$$
P_\gamma(f)(z)= \int_{\B_n} f(w)K_\gamma(z, w) dv(w).
$$
 Fix  $\gamma$ so that $\gamma>\alpha$ and $2(\gamma+1)>\beta+1$, it follows from \cite[Theorem 2.11]{Zh05} that $P_\gamma$ is bounded from $L^1_\alpha$ to $A^1_\alpha$ and bounded from $L^2_\beta$ to $A^2_\beta$. Set $U={\rm Id}_4\circ {\rm Id}_2\circ P_\gamma$, then $U$ is bounded from $L_\alpha^{1}$ to ${L}_\beta^{2}$ with $\|U\|_{L_\alpha^1 \rightarrow  {L}_\beta^{2}}\lesssim 1$. By \cite[Theorem 3.4]{DJT95} we get $U: {L}_\alpha^{1} \rightarrow  {L}_\beta^{2}$ is 1-summing with $\pi_1(U)\lesssim \|U\|_{ L_\alpha^1 \rightarrow  L_\beta^2}$. As mentioned on page 37 in \cite{DJT95},  the restriction map $U|_{A_\alpha^{1}}$ is also 1-summing with
$$
\pi_1(U|_{A_\alpha^{1}}: A_\alpha^1 \rightarrow  L_\beta^2)\le \pi_1(U)\lesssim 1.
$$
Therefore,
$$
{\rm Id_2}(A_\alpha^1 \rightarrow  A_\beta ^2)=P_\gamma({L}_\beta^{2} \rightarrow  A_\beta^{2})\circ U|_{{A}_\alpha^{1}}(A_\alpha^{1} \rightarrow  {L}_{\beta}^{2})
$$
 is 1-summing with
 \begin{equation}\label{fp-d}
 \pi_1({\rm Id_2}: A_\alpha^1 \rightarrow  A_\beta^2)\lesssim 1.
\end{equation}

 It is easy to see that $T_\mu= R_{0, N}\circ {\rm Id_3} \circ {\rm Id_2} \circ R^{0, N} \circ T_\mu $. Then we know is $T_\mu: A^p \to A^p $ is $r$-summing. Furthermore, by (\ref{1<r<p'}), (\ref{fp-a}), (\ref{fp-b}), (\ref{fp-d}) and the ideal property we get
$$
\pi_r(T_\mu)\le \pi_1(T_\mu)\lesssim\|\widetilde{\mu}\|_{L^{p'}(d\lambda)}.
$$
This and (\ref{4.3333}) give the norm equivalence (\ref{norm4-7}).
 \end{proof}

 \subsection{The $r$-summing Toeplitz operators on $A^{p}$ for $p\ge 2$ and $p' \leq  r < p$.} \label{4.3}

\begin{proposition}\label{prop3.3}
 Suppose  $p\ge 2$ and $p' \leq  r < p$. Let $\mu$ be a positive Borel measure on \(\B_n\). Then  \({T}_{\mu} : A^{p} \rightarrow  A^{p}\) is \(r\)-summing if and only if  \(\widetilde{\mu} \in L^ {r}(d\lambda)\). Furthermore,
\begin{equation}\label{norm1}
\pi_r(T_\mu)\simeq\|\widetilde{\mu}\|_{L^r(d\lambda)}.
 \end{equation}
\end{proposition}

\begin{proof}
Suppose \({T}_{\mu} \in \Pi_r( A^{p}, A^{p})\). Set
$$
\Phi: \B_n  \longrightarrow  A^{p},
$$
$$
\ \ \ \ \ \ \ \ \ z   \longmapsto  k_{z, p}(\cdot)
$$
Lemma \ref{lem8} gives
$$
\int_{\B_n}\|T_\mu  (k_{z, p})\|_{A^{p}}^{r}{d\lambda(z)} \leq  \left({\pi }_{r}{\left( T_\mu\right) }\right)^{r}\cdot\mathop{\sup }\limits_{{g  \in  {B}_{A^{p'}}}}\int_{\B_n}{\left|\langle g, k_{z, p} \rangle\right| }^{r}d\lambda(z).
$$
On one hand, the sub-mean value property of $|T_\mu (k_{z, p})|^p$ yields
$$
\widetilde{\mu}(z)^p\lesssim\|T_\mu (k_{z,p})\|_{A^{p}}^{p},
$$
which follows
$$
\int_{\B_n}\widetilde{\mu}(z)^{r}{d\lambda(z)}\lesssim \int_{\B_n}\|T_\mu  (k_{z, p})\|_{A^{p}}^{r}{d\lambda(z)}.
$$
On the other hand, since $r\ge p'$, (\ref{formula1}) gives
\beqm
\mathop{\sup }\limits_{g  \in  B_{A^{p'}}}\int_{\B_n}{\left|\langle g, k_{z, p} \rangle\right| }^{r}{d\lambda(z)}
&\simeq& \mathop{\sup }\limits_{{g  \in  {B}_{A^{p'}}}}\int_{\B_n}{\left|g(z)\cdot K(z,z)^{\fr 1p-1}\right| }^{r}d\lambda(z) \\
&\lesssim&  \mathop{\sup }\limits_{{g  \in  B_{A^{p'}}}}\int_{\B_n}|g(z)|^{p'}\cdot\|g\|^{r-{p'}}_{A^{p'}}K(z,z)^{\fr {r-p'}{p'}}\cdot K(z,z)^{\fr rp-r+1}dv(z)\\
&=& \mathop{\sup }\limits_{{g  \in  B_{A^{p'}}}}\|g\|^{r}_{A^{p'}}\\
&= &1.
\eqm
The above estimates show that
 \begin{equation}\label{mut-bounded}
\int_{\B_n}\widetilde{\mu}(z)^{r}{d\lambda(z)}\lesssim \left({\pi }_{r}{( T_\mu)}\right)^{r},
 \end{equation}
which tells us \(\widetilde{\mu}\in  {L}^{r}(d\lambda)\) with $\|\widetilde{\mu}\|_{L^r(d\lambda)}\lesssim \pi_r(T_\mu)$.

Conversely, suppose \(\widetilde{\mu} \in  {L}^{r}(d\lambda)\). Similar to the proof of the previous subsection, take $N$ to be positive so that
$$
  \eta= Nr-(n+1)\left (1-\frac rp\right)>0.
$$
Consider the chain of maps:
$$
%A^p\stackrel{{\rm Id}_1}{\longrightarrow} A^\infty_{(n+)/p}
A^p\stackrel{R^{0, N}\circ T_\mu}{\longrightarrow}A_\eta^r\stackrel{{\rm Id}}{\longrightarrow}{A^p_{Np}}\stackrel{R_{0, N}}{\longrightarrow}A^p,
$$
We claim that the following two statements hold.

(i)  ${\rm Id}:  A^r_\eta \to A^p_{Np}$ and $R_{0, N}: A^p_{Np}\to A^p$ are bounded.

(ii)  $R^{0, N}\circ T_\mu: A^p\to L_\eta^r$ is order bounded, and then $R^{0, N}\circ T_\mu: A^p\to A_\eta^r$ is $r$-summing.\\

In fact,  for $f\in A^r_\eta $, by (\ref{weighted-a}) we know
 $$
   |f(z)| \le  \|f\|_{A^r_\eta} (1-|z|^2)^{-\frac {n+1+\eta} r}.
 $$
Notice that $ Np-\frac {n+1+\eta} r (p-r)= \eta$,     then
\beqm
  \|f\|_{A^p_{Np}}^p
   &\le& %\|f\|_{A^r_\beta}
   \int_{\B_n} \left (|f(w)| (1-|w|^2)^ {\frac {n+1+\eta} r}\right)^{p-r} \left( |f(w)|^r (1-|w|^2)^{Np-\frac {n+1+\eta} r (p-r)}\right) dv(w)\\
    &\lesssim &\|f\|_{A^r_\eta}^{p-r} \int_{\B_n} |f(w)|^r (1-|w|^2)^{\eta} dv(w).
\eqm
With this we have ${\rm Id}:  A^r_\eta \to A^p_{Np}$ is bounded
$$
 \|{\rm Id}\|_{ A^r_\eta \to A^p_{Np}}\lesssim 1.
$$
The boundedness of $R_{0, N}: A^p_{Np}\to A^p$ comes from (\ref{fp-c}).

Now we prove that  $R^{0, N}\circ T_\mu: A^p\to L_\eta^r$ is order bounded.
For $f\in B_{A^p}$, similar to that in the previous subsection we have
\beqm
     & &\left | R^{0, N}\circ T_\mu f(z)\right|\\
 &\lesssim &
      \int_{\B_n}\frac{|f(w)|\widehat{\mu}_\delta(w) }{|1-\langle z,w\rangle|^{n+1 +N}} dv(w)\\
 &\lesssim & \|f\|_{A^p} \int_{\B_n} \widehat{\mu}_\delta(w) K(w,
     w)^{\frac 1 r}  \frac {(1-|w|^2)^{(n+1) \left(\frac 1 r- \frac 1 p\right)} }  {|1-\langle z,w\rangle|^{n+1+N}} dv(w)\\
 &=&\|f\|_{A^p} \cdot S\left( \widehat{\mu}_\delta K(\cdot,
     \cdot)^{\frac 1 r} \right)(z),
\eqm
where $S=S_{b,c}$ with $b=(n+1) \left(\frac 1 r- \frac 1 p\right)$ and $ c=n+1+N$ is  the operator defined in (\ref{S_bc}).
Applying \cite[Theorem 3]{Zhao15} we know $S$ is bounded from $L^r$ to $L^r_{\eta }$ with the estimate
  $$
    \left \| S\left( \widehat{\mu}_\delta K(\cdot,
     \cdot)^{\frac 1 r} \right) \right \|_{L^r_{\eta }} \le \|\widehat{\mu}_\delta K(\cdot,
     \cdot)^{\frac 1 r}\|_{L^r}\cdot\|S\|  = \| \widehat{\mu}_\delta\|_{L^r(d \lambda)}\cdot\|S\|.
  $$
Hence,  $R^{0, N}\circ T_\mu:  A^p\to L_\eta^r $ is order bounded. And
$$
 {\left\|\mathop{\sup }\limits_{f \in  B_{A^p}}\left |R^{0, N}\circ T_\mu  f\right | \right\|}_{L_\eta^r} \lesssim  \| \widehat{\mu}_\delta\|_{L^r(d \lambda)}.
$$
Lemma \ref {lem66} tells us $R^{0, N}\circ T_\mu:  A^p\to L_\eta^r $ is $r$-summing with the estimate
$$%\begin{equation}\label{id-b}
  \pi_r\left (R^{0, N}\circ T_\mu: A^p\to L_\eta^r \right)  \lesssim   {\left\|\mathop{\sup }\limits_{f \in  B_{A^p}}\left |R^{0, N}\circ T_\mu  f\right | \right\|}_{L_\eta^r} \lesssim  \| \widehat{\mu}_\delta\|_{L^r(d \lambda)}.
$$%\end{equation}.
Since  $R^{0, N}\circ T_\mu(A^p)\subset A_\eta^r$,  we obtain
$$
  \pi_r\left (R^{0, N}\circ T_\mu: A^p\to A_\eta^r \right)  \lesssim   \| \widehat{\mu}_\delta\|_{L^r(d \lambda)}.
$$
Therefore, $T_\mu = R_{0, N}\circ {\rm Id}  \circ R^{0, N} \circ T_\mu $ is $r$-summing on $A^p$. Moreover, 
\begin{equation}\label{id-b}
  \pi_r\left (T_\mu: A^p\to A^p \right)  \lesssim   \| \widehat{\mu}_\delta\|_{L^r(d \lambda)}.
\end{equation}

 The norm equivalence (\ref{norm1}) comes from (\ref{mut-bounded}) and (\ref{id-b}).
 \end{proof}

\subsection{The $r$-summing Toeplitz operators on $A^{p}$ for $p\ge 2$ and $r \geq p$.}

\begin{proposition}\label{prop3.2}
Suppose $p\ge 2$ and $r \geq p$. Let $\mu$ be a positive Borel measure on \(\B_n\). Then  \({T}_{\mu} : A^{p} \rightarrow  A^{p}\) is \(r\)-summing if and only if  \(\widetilde{\mu} \in L^ {p}(d\lambda)\). Furthermore,
\begin{equation}\label{norm3319}
\pi_r(T_\mu)\simeq\|\widetilde{\mu}\|_{L^p(d\lambda)}.
 \end{equation}
\end{proposition}

\begin{proof}
Suppose \(\widetilde{\mu} \in  {L}^{p}(d\lambda)\), we claim $T_\mu: A^p\to  L^p$ is order bounded. To see this, for $\delta>0$ fixed, we set
$$
h(z)=\int_{\B_n}K(\xi,\xi)^{\fr1p}\widehat{\mu}_\delta(\xi)
|K(z,\xi)|dv(\xi).
$$
The boundedness of the  operator $S_{0, n+1}$ defined in (\ref{S_bc}) implies that
 \begin{equation}\label{intest}
\int_{\B_n}h(z)^{p}dv(z)\lesssim \int_{\B_n}K(z,z)\widehat{\mu}_\delta(z)^pdv(z)= \|\widehat{\mu}_\delta\|_{L^p(d\lambda)}^p <\infty.
\end{equation}
Now for $g\in  A^p$, (\ref{formula1}) gives
\beqm
{\left| T_\mu(g)(z) \right| }
&\lesssim& \int_{\B_n}|g(\xi)||K(z,\xi)|\widehat{\mu}_\delta(\xi)dv(\xi)\\
&\lesssim & \|g\|_{A^p}\cdot\int_{\B_n}K(\xi,\xi)^{\fr1p}\widehat{\mu}
_\delta(\xi)|K(z,\xi)|dv(\xi)\\
&=& \|g\|_{A^p} h(z) .
\eqm
Lemma \ref{lem66}, together with (\ref{cotype}) and (\ref{intest}), gives $T_\mu\in \Pi_p( A^p, L^p)$. Since $T_\mu(A^p)\subset A^p$, we  have $T_\mu\in \Pi_p( A^p, A^p)$,
  \begin{equation}\label{3110}
\pi_p(T_\mu: A^p\to A^p)= \pi_p(T_\mu: A^p\to L^p)\le  {\left\|\mathop{\sup }\limits_{g \in  B_{A^p}}|T_\mu(g)| \right\|}_{{L}^{p}}\lesssim  \|\widetilde{\mu}\|_{L^p(d\lambda)}.
\end{equation}
%Hence  $T_\mu\in \Pi_r( A^{p}, A^p)$ since $r\ge p$. Furthermore,

%\pi_r(T_\mu)\le\pi_p(T_\mu) \lesssim \|\widetilde{\mu}\|_{L^p(d\lambda)}.

Conversely, suppose $T_\mu\in \Pi_r( A^{p}, A^p)$. Since the cotype of $A^{p}$ is $p \geq  2$, it follows from \cite[Theorem 11.13]{DJT95} that $T_\mu\in \Pi_{p,2}(A^p, A^p)$. Precisely, for $\{f_j\}_{j=1}^m\subset A^{p}$, there holds
\beqm
\left(\sum_{j=1}^{m}\|T_\mu f_j\|_{A^{p}}^{p}\right)^{\fr 1 p}
&\lesssim& \pi_r(T_\mu) \cdot \mathop{\sup }\limits_{{g  \in  {B}_{A^{p'}}}}\left(\sum_{j=1}^{m}\left|\langle g, f_j \rangle\right|^2\right)^{\fr 1 2}\\
&\lesssim &  \pi_r(T_\mu)\cdot \mathop{\sup }\limits_{{g  \in  {B}_{A^{p'}}}}\left(\sum_{j=1}^{m}\left|\langle g, f_j \rangle\right|^{p'}\right)^{\fr 1 {p'}}.
\eqm
Fix some  $\delta$-lattice $\{a_j\}$ with $\delta$ small, set $f_j=k_{a_j, p}$, then
$\widetilde{\mu}(z_j)\lesssim \|T_\mu (k_{a_j, p})\|_{A^{p}}$.
From the above we obtain
\beqm
\left(\sum_{j=1}^{m}\widetilde{\mu}(a_j)^{p}\right)^{\fr 1 p}
&\lesssim& \pi_r(T_\mu)\cdot \mathop{\sup }\limits_{{g  \in  {B}_{A^{p'}}}}\left(\sum_{j=1}^{m}\left|\langle g, k_{a_j, p} \rangle\right|^{p'}\right)^{\fr 1 {p'}}\\
&\simeq&\pi_r(T_\mu)\cdot \mathop{\sup }\limits_{{g  \in  {B}_{A^{p'}}}}\left(\sum_{j=1}^{m}\left|g(a_j)K(a_j,a_j)^{1-\fr 1 p} \right|^{p'}\right)^{\fr 1 {p'}}\\
&\lesssim& \pi_r(T_\mu)\cdot \mathop{\sup }\limits_{{g  \in  {B}_{A^{p'}}}}\|g\|_{A^{p'}}^{p'}\\
&= & \pi_r(T_\mu).
\eqm
Since $\widehat{\mu}_\delta(a_j)\le C \widetilde{\mu}(a_j)$,  Lemma \ref{lem2.3} tells us that  \(\widetilde{\mu} \in  {L}^{p}(d\lambda)\) with  $$
\|\widetilde{\mu}\|_{L^{p}(d\lambda)}\simeq\|\{\widehat{\mu}_\delta(a_j)\}\|_{l^{p}}\lesssim \pi_r(T_\mu).
$$

This along with (\ref{3110}) gives the norm equivalence (\ref{norm3319}).
 \end{proof}

Summarizing the above 4 subsections, we have obtain  the proof of Theorem \ref{THM2} for all $1< p<\infty$ and $r\ge 1.$

\section{Further Remarks}\label{remark}

Take $1\le p=q\le 2$, and take $r=1$ in  Theorem \ref{THM1}.  We know that the Carleson embedding $J_\mu: A^p\to L^p(\mu)$ is absolutely summing if and only if
$$
\widetilde{\mu}\in L^1(d\lambda).
$$
In other words,  the absolutely summing for $J_\mu$ on $A^p$ can be formulated with the same condition.
For the Toeplitz operators, Theorem \ref{THM2} says  that $T_\mu$
is absolutely summing on $A^p$, $1<p\le 2$,  if and only if
\begin{equation}\label{Rem-1}
\widetilde{\mu}\in L^2 (d\lambda).
\end{equation}
Comparing these observations with the known characterizations of $T_\mu$ and $J_\mu$ bounded (or compact)  on $A^p$, we have good reason to believe that the condition (\ref{Rem-1}) is still far from being necessary and sufficient for $T_\mu$  to be absolutely summing on $A^1$. It is therefore natural to ask what conditions characterize $T_\mu\in \pi_1(A^1\to A^1)$. At present we do not even have a plausible conjecture.\\

The approach in this paper works well on weighted Bergman spaces. Precisely, for $\alpha>-1$ let $K_\alpha(\cdot, \cdot)$ be the Bergman kernel on $A^2_\alpha$. Similar to  (\ref{toep}),   we define the Toeplitz operator $T_{\mu, \alpha}$ on $A^p_\alpha$ as
 $$
    T_{\mu, \alpha}  f(z) = \int_{\B_n} f(w)K_\alpha (z,w)d\mu(w).  
 $$
  With these adjustment we have the following theorem.

\begin{theorem}\label{THM3}
 Let $\mu$ be a positive Borel measure on $\B$. For $1< p<\infty$ and $ 1\le r<\infty$,
the Toeplitz operator   \(T_{\mu, \alpha}: A^p_\alpha \rightarrow  A^p_\alpha \) is \(r\)-summing if and only if \(\widetilde{\mu} \in  {L}^{\kappa}(d\lambda)\). Furthermore,
$$
\pi_r(T_\mu : A^p_\alpha \rightarrow  A^p_\alpha ) \simeq  \|\widetilde{\mu}\|_{L^{\kappa}(d\lambda)}.
$$
\end{theorem}

 How should absolutely summing positive Toeplitz operators be characterized on Bergman spaces of more general domains $\Omega\subset \mathbb{C}^n$?
On bounded symmetric domains the entire machinery of this paper transfers verbatim: the transitive group action preserves all kernel estimates, and the fractional-derivative and fractional-integral operators $R^{0,N}, R_{0,N}$ are available in the same way as that on $\B_n$. For bounded strongly pseudoconvex domains $\Omega$ the $r$-summing Carleson embedding $J_\mu:A^p(\Omega)\to L^p(\Omega)$ can still be established by the same strategy, and the Bergman-kernel estimates we need are already known. What is missing, however, is the exact fractional-derivative/fractional-integral operators used in Subsections \ref{4.2} and \ref{4.3}. To prove the analogue of Theorem \ref{THM2} we expect that suitable substitutes-most likely built on higher-order gradients-can be constructed to close the gap.

\begin{center}

\noindent
\end{center}

\end{document}